\documentclass{amsart}%
\usepackage{amsfonts}
\usepackage{graphicx}
\usepackage{amscd}%
\usepackage{amsmath}%
\usepackage{amssymb}
\newtheorem{theorem}{Theorem}
\theoremstyle{plain}

\newtheorem{case}{Case}

\newtheorem{conjecture}{Conjecture}
\newtheorem{corollary}{Corollary}

\newtheorem{lemma}{Lemma}

\newtheorem{proposition}{Proposition}
\newtheorem{remark}{Remark}

\numberwithin{equation}{section}
\numberwithin{theorem}{section}
\numberwithin{algorithm}{section}
\numberwithin{axiom}{section}
\numberwithin{case}{section}
\numberwithin{claim}{section}
\numberwithin{conclusion}{section}
\numberwithin{condition}{section}
\numberwithin{conjecture}{section}
\numberwithin{corollary}{section}
\numberwithin{criterion}{section}
\numberwithin{definition}{section}
\numberwithin{example}{section}
\numberwithin{exercise}{section}
\numberwithin{lemma}{section}
\numberwithin{notation}{section}
\numberwithin{problem}{section}
\numberwithin{proposition}{section}
\numberwithin{remark}{section}
\numberwithin{solution}{section}

\begin{document}
\title{An integral equation in conformal geometry}
\author{Fengbo Hang}
\address{Department of Mathematics, Princeton University, Fine Hall, Washington Road,
Princeton, NJ\ 08544}
\email{fhang@math.princeton.edu}
\author{Xiaodong Wang}
\address{Department of Mathematics, Michigan State University, East Lansing, MI 48824}
\email{xwang@math.msu.edu}
\author{Xiaodong Yan}
\address{Department of Mathematics, Michigan State University, East Lansing, MI 48824}
\email{xiayan@math.msu.edu}
\maketitle

\section{Introduction\label{sec1}}

Among the many proofs of two dimensional isoperimetric inequalities, the one
due to Carleman \cite{C} is particularly interesting. Indeed by an application
of Riemann mapping theorem we only need to show%
\begin{equation}
\int_{D}e^{2u}dx\leq\frac{1}{4\pi}\left(  \int_{S^{1}}e^{u}d\theta\right)
^{2}\label{eq1.1}%
\end{equation}
for every harmonic function $u$ on $D$. Here $D$ is the unit disk in the
plane. Carleman deduced (\ref{eq1.1}) by showing%
\[
\int_{D}\left\vert f\right\vert ^{2}dx\leq\frac{1}{4\pi}\left(  \int_{S^{1}%
}\left\vert f\right\vert d\theta\right)  ^{2}%
\]
for every holomorphic function $f$ on $D$. Along this line, in \cite{J} Jacobs
showed that for every bounded open subset $\Omega$ of $\mathbb{R}^{2} $ with
smooth boundary, there exists a positive constant $c_{\Omega}$ such that for
every holomorphic function $f$ on $\Omega$,%
\[
\int_{\Omega}\left\vert f\right\vert ^{2}dx\leq c_{\Omega}\left(
\int_{\partial\Omega}\left\vert f\right\vert ds\right)  ^{2}.
\]
Moreover when $\Omega$ is not simply connected, the best constant $c_{\Omega
}>\frac{1}{4\pi}$ and it is achieved by some particular holomorphic function
$f$. Here we formulate a higher dimensional generalization of these statements.

Assume $n\geq3$, $\left(  M^{n},g\right)  $ is a smooth compact Riemannian
manifold with nonempty boundary $\Sigma=\partial M$, we write the
isoperimetric ratio%
\begin{equation}
I\left(  M,g\right)  =\frac{\left\vert M\right\vert ^{\frac{1}{n}}}{\left\vert
\Sigma\right\vert ^{\frac{1}{n-1}}}.\label{eq1.2}%
\end{equation}
Here $\left\vert M\right\vert $ is the volume of $M$ with respect to $g$ and
$\left\vert \Sigma\right\vert $ is the area of $\Sigma$. Let $\left[
g\right]  =\left\{  \rho^{2}g:\rho\in C^{\infty}\left(  M\right)
,\rho>0\right\}  $ be the conformal class of $g$. The set%
\[
\left\{  \widetilde{g}\in\left[  g\right]  :\text{the scalar curvature
}\widetilde{R}=0\right\}
\]
is nonempty if and only if the first eigenvalue of the conformal Laplacian
operator $L_{g}=-\frac{4\left(  n-1\right)  }{n-2}\Delta+R$ with respect to
Dirichlet boundary condition, $\lambda_{1}\left(  L_{g}\right)  $ is strictly
positive (see Section \ref{sec2}).

Assume $\lambda_{1}\left(  L_{g}\right)  >0$, we denote%
\begin{equation}
\Theta_{M,g}=\sup\left\{  I\left(  M,\widetilde{g}\right)  :\widetilde{g}%
\in\left[  g\right]  \text{ with }\widetilde{R}=0\right\}  .\label{eq1.3}%
\end{equation}
Standard technique from harmonic analysis gives us $\Theta_{M,g}<\infty$ (see
Proposition \ref{prop2.1}). But is $\Theta_{M,g}$ achieved? In another word,
can we find a conformal metric with zero scalar curvature maximizing the
isoperimetric ratio?

It follows from \cite[theorem 1.1]{HWY} or Theorem \ref{thm3.1} that%
\[
\Theta_{\overline{B}_{1},g_{\mathbb{R}^{n}}}=I\left(  \overline{B}%
_{1},g_{\mathbb{R}^{n}}\right)  =n^{-\frac{1}{n-1}}\omega_{n}^{-\frac
{1}{n\left(  n-1\right)  }},
\]
here $\omega_{n}$ is the volume of the unit ball in $\mathbb{R}^{n}$ and
$g_{\mathbb{R}^{n}}$ is the Euclidean metric on $\mathbb{R}^{n}$. This just
says that $\Theta_{\overline{B}_{1},g_{\mathbb{R}^{n}}}$ is achieved by the
standard metric. In general we have the following

\begin{theorem}
\label{thm1.1}Assume $n\geq3$, $\left(  M^{n},g\right)  $ is a smooth compact
Riemannian manifold with nonempty boundary and $\lambda_{1}\left(
L_{g}\right)  >0$, then%
\[
n^{-\frac{1}{n-1}}\omega_{n}^{-\frac{1}{n\left(  n-1\right)  }}=\Theta
_{\overline{B}_{1},g_{\mathbb{R}^{n}}}\leq\Theta_{M,g}<\infty.
\]
If in addition $\Theta_{\overline{B}_{1},g_{\mathbb{R}^{n}}}<\Theta_{M,g}$,
then $\Theta_{M,g}$ is achieved by some conformal metrics with zero scalar curvature.
\end{theorem}

The problem illustrates very similar behavior as the Yamabe problem of finding
constant scalar curvature metrics in a fixed conformal class (cf. \cite{LP})
and its boundary versions (cf. \cite{E1,E2}). On the other hand, it has more
nonlocal features (e.g. the Euler-Lagrange equation is a nonlinear integral
equation) than the two well studied problems. In analogy with the solution of
the Yamabe problem, we make the following conjecture.

\begin{conjecture}
\label{conj1.1}Assume $n\geq3$, $\left(  M^{n},g\right)  $ is a smooth compact
Riemannian manifold with nonempty boundary and $\lambda_{1}\left(
L_{g}\right)  >0$. If $\left(  M,g\right)  $ is not conformally diffeomorphic
to $\left(  \overline{B}_{1},g_{\mathbb{R}^{n}}\right)  $, then $\Theta
_{M,g}>\Theta_{\overline{B}_{1},g_{\mathbb{R}^{n}}}$.
\end{conjecture}

In Section \ref{sec2} below, we will describe some basics related to the above
problem and reformulate it as a maximization problem for harmonic extensions.
We will also discuss some elementary estimates of the Poisson kernels and show
$\Theta_{M,g}$ is always finite. In Section \ref{sec3} we will show
$\Theta_{\overline{B}_{1},g_{\mathbb{R}^{n}}}$ is achieved by the standard
metric itself and deduce some corollaries. This is a consequence of
\cite[theorem 1.1]{HWY}. However the approach we present here is different and
of independent interest. In Section \ref{sec4} we will prove the regularity of
the solutions to the Euler-Lagrange equations of the maximization problem for
harmonic extensions. In Section \ref{sec5} we derive some asymptotic expansion
formulas for the standard Poisson kernel and the Poisson kernel for the
conformal Laplacian operators. These expansion formulas will be useful in the
future study of Conjecture \ref{conj1.1}. In Section \ref{sec6}, we will
derive the concentration compactness principle for the maximization problem
and this will be used in the last section to deduce Theorem \ref{thm1.1}.

\textbf{Acknowledgment}: The research of F. Hang is supported by National
Science Foundation Grant DMS-0647010 and a Sloan Research Fellowship. The
research of X. Wang is supported by National Science Foundation Grant
DMS-0505645. The research of X. Yan is supported by National Science
Foundation Grant DMS-0401048. We would like to thank Professor R. Mazzeo and
Professor E. Stein for some helpful discussions. We also thank the referee for
his/her careful reading of the original manuscript.

\section{Some preparations\label{sec2}}

Assume $n\geq3$, $\left(  M^{n},g\right)  $ is a smooth compact Riemannian
manifold with boundary $\Sigma=\partial M$. The conformal Laplacian operator
is given by%
\[
L_{g}=-\frac{4\left(  n-1\right)  }{n-2}\Delta+R.
\]
It satisfies the transformation law%
\[
L_{\rho^{\frac{4}{n-2}}g}\varphi=\rho^{-\frac{n+2}{n-2}}L_{g}\left(
\rho\varphi\right)  \quad\text{for }\rho,\varphi\in C^{\infty}\left(
M\right)  ,\rho>0\text{.}%
\]
Let%
\[
E_{g}\left(  \varphi,\psi\right)  =\int_{M}\left[  \frac{4\left(  n-1\right)
}{n-2}\nabla\varphi\cdot\nabla\psi+R\varphi\psi\right]  d\mu,\quad
E_{g}\left(  \varphi\right)  =E_{g}\left(  \varphi,\varphi\right)  ,
\]
here $d\mu$ is the measure generated by $g$, then it follows from the
transformation law that%
\begin{equation}
E_{\rho^{\frac{4}{n-2}}g}\left(  \varphi\right)  =E_{g}\left(  \rho
\varphi\right)  \quad\text{for }\rho,\varphi\in C^{\infty}\left(  M\right)
,\rho>0,\left.  \varphi\right\vert _{\Sigma}=0.\label{eq2.1}%
\end{equation}
Let $\lambda_{1}\left(  L_{g}\right)  $ be the first eigenvalue of $L_{g}$
with respect to the Dirichlet boundary condition, then%
\[
\lambda_{1}\left(  L_{g}\right)  =\inf_{\varphi\in H_{0}^{1}\left(  M\right)
\backslash\left\{  0\right\}  }\frac{E_{g}\left(  \varphi\right)  }{\int
_{M}\varphi^{2}d\mu_{g}}.
\]
Assume $\rho\in C^{\infty}\left(  M\right)  ,$ $\rho>0$. It follows from
(\ref{eq2.1}) that $\lambda_{1}\left(  L_{g}\right)  <0$ implies $\lambda
_{1}\left(  L_{\rho^{\frac{4}{n-2}}g}\right)  <0$. On the other hand, if
$\lambda_{1}\left(  L_{g}\right)  \geq0$, then $\lambda_{1}\left(
L_{\rho^{\frac{4}{n-2}}g}\right)  \geq\left(  \max_{M}\rho\right)  ^{-\frac
{4}{n-2}}\lambda_{1}\left(  L_{g}\right)  $. Hence the sign of the first
eigenvalue of the conformal Laplacian operator does not depend on the choice
of particular metric in a conformal class. This sign is useful because of the
following fact: $\lambda_{1}\left(  L_{g}\right)  >0$ if and only if we may
find a scalar flat metric in the conformal class of $g$. The only thing we
need to verify is we can find a scalar flat conformal metric when $\lambda
_{1}\left(  L_{g}\right)  >0$. To see this we may solve the Dirichlet problem%
\[
\left\{
\begin{array}
[c]{l}%
L_{g}\rho=0\text{ on }M\\
\left.  \rho\right\vert _{\Sigma}=1
\end{array}
\right.  .
\]
We claim $\rho>0$ on $M$. To see this, we let $\varphi$ be the first
eigenfunction of $L_{g}$ with $\varphi>0$ on $M\backslash\Sigma$ and $\left.
\varphi\right\vert _{\Sigma}=0$. Let $w=\frac{\rho}{\varphi}$, then%
\[
-\frac{4\left(  n-1\right)  }{n-2}\Delta w-\frac{8\left(  n-1\right)  }%
{n-2}\frac{\nabla\varphi}{\varphi}\cdot\nabla w+\lambda_{1}w=0\text{ on
}M\backslash\Sigma.
\]
Since $w\left(  x\right)  \rightarrow\infty$ as $x\rightarrow\Sigma$, it
follows from strong maximum principle that $w>0$ on $M\backslash\Sigma$, hence
$\rho>0$ on $M$. Note that $R_{\rho^{\frac{4}{n-2}}g}=\rho^{-\frac{n+2}{n-2}%
}L_{g}\rho=0$, we find the needed metric.

Assume $\lambda_{1}\left(  L_{g}\right)  >0$, the Green's function $G_{L}$ of
$L_{g}$ satisfies%
\[
\left\{
\begin{array}
[c]{l}%
\left(  L_{g}\right)  _{x}G_{L}\left(  x,y\right)  =\delta_{y}\text{ on }M,\\
G_{L}\left(  x,y\right)  =0\text{ for }x\in\Sigma.
\end{array}
\right.
\]
The Poisson kernel of $L_{g}$ is given by%
\[
P_{L}\left(  x,\xi\right)  =-\left.  \frac{4\left(  n-1\right)  }{n-2}%
\frac{\partial G_{L}\left(  x,y\right)  }{\partial_{y}\nu}\right\vert _{y=\xi
},
\]
here $\nu$ is the unit outer normal direction. The solution of $\left\{
\begin{array}
[c]{l}%
L_{g}u=0\text{ on }M\\
\left.  u\right\vert _{\Sigma}=f
\end{array}
\right.  $ is given by $u\left(  x\right)  =\left(  P_{L}f\right)  \left(
x\right)  =\int_{\Sigma}P_{L}\left(  x,\xi\right)  f\left(  \xi\right)
dS\left(  \xi\right)  $, here $dS$ is the measure generated by $g$ on $\Sigma
$. If $\rho$ is a positive smooth function, then we have the following
transformation laws,%
\[
G_{L,\rho^{\frac{4}{n-2}}g}\left(  x,y\right)  =\frac{G_{L,g}\left(
x,y\right)  }{\rho\left(  x\right)  \rho\left(  y\right)  },\quad
P_{L,\rho^{\frac{4}{n-2}}g}\left(  x,\xi\right)  =\frac{P_{L,g}\left(
x,\xi\right)  }{\rho\left(  x\right)  \rho\left(  \xi\right)  ^{\frac{n}{n-2}%
}},
\]
and%
\[
P_{L,\rho^{\frac{4}{n-2}}g}f=\rho^{-1}P_{L,g}\left(  \rho f\right)  .
\]
If $\widetilde{g}\in\left[  g\right]  $ has zero scalar curvature, then
$\widetilde{g}=u^{\frac{4}{n-2}}g$ for some positive smooth function $u$ on $M
$ with $L_{g}u=0$. Let $f=\left.  u\right\vert _{\Sigma}$, then $u=P_{L}f$ and%
\[
I\left(  M,\widetilde{g}\right)  =\frac{\left\vert P_{L}f\right\vert
_{L^{\frac{2n}{n-2}}\left(  M\right)  }^{\frac{2}{n-2}}}{\left\vert
f\right\vert _{L^{\frac{2\left(  n-1\right)  }{n-2}}\left(  \Sigma\right)
}^{\frac{2}{n-2}}}.
\]
Hence%
\begin{align}
& \Theta_{M,g}\label{eq2.2}\\
& =\sup\left\{  \frac{\left\vert P_{L}f\right\vert _{L^{\frac{2n}{n-2}}\left(
M\right)  }^{\frac{2}{n-2}}}{\left\vert f\right\vert _{L^{\frac{2\left(
n-1\right)  }{n-2}}\left(  \Sigma\right)  }^{\frac{2}{n-2}}}:f\in C^{\infty
}\left(  \Sigma\right)  ,f>0\right\} \nonumber\\
& =\sup\left\{  \frac{\left\vert P_{L}f\right\vert _{L^{\frac{2n}{n-2}}\left(
M\right)  }^{\frac{2}{n-2}}}{\left\vert f\right\vert _{L^{\frac{2\left(
n-1\right)  }{n-2}}\left(  \Sigma\right)  }^{\frac{2}{n-2}}}:f\in
L^{\frac{2\left(  n-1\right)  }{n-2}}\left(  \Sigma\right)  ,f\neq0\right\}
\nonumber\\
& =\left[  \sup\left\{  \left\vert P_{L}f\right\vert _{L^{\frac{2n}{n-2}%
}\left(  M\right)  }:f\in L^{\frac{2\left(  n-1\right)  }{n-2}}\left(
\Sigma\right)  ,\left\vert f\right\vert _{L^{\frac{2\left(  n-1\right)  }%
{n-2}}\left(  \Sigma\right)  }=1\right\}  \right]  ^{\frac{2}{n-2}}.\nonumber
\end{align}
The second equality above follows from the fact $P_{L}$ is positive and an
approximation procedure.

It follows easily from the definition of $\Theta_{M,g}$ (see (\ref{eq1.3}))
that $\Theta_{M,g}$ depends only on $\left[  g\right]  $. As a consequence we
may choose the background metric $g$ with zero scalar curvature. Under this
assumption the conformal Laplacian operator reduces to the constant multiple
of the Laplacian operator. To continue we will need some estimates of the
Poisson kernels.

\subsection{Basic estimates for Poisson kernel and harmonic
extensions\label{sec2.1}}

Let us fix some notations. Throughout this subsection, we always assume
$n\geq2$, $\left(  M^{n},g\right)  $ is a smooth compact Riemannian manifold
with boundary $\Sigma=\partial M$. For convenience we fix a smooth compact
Riemannian manifold without boundary, $\left(  \overline{M}^{n},g\right)  $
such that $\left(  M,g\right)  $ is a smooth domain in $\left(  \overline
{M},g\right)  $. Denote $d$ as the distance on $\overline{M}$ generated by $g$
and $d_{\Sigma}$ as the distance on $\Sigma$ (when $\Sigma$ is not connected
and $\xi_{1},\xi_{2}\in\Sigma$ lie in different components, we set $d_{\Sigma
}\left(  \xi_{1},\xi_{2}\right)  $ equal to the maximal diameter of all the
components of $\Sigma$). We write $t=t\left(  x\right)  =d\left(
x,\Sigma\right)  $ for $x\in\overline{M}$. Assume $\delta_{0}>0$ is small
enough such that $V=\left\{  x\in\overline{M}:t\left(  x\right)  <2\delta
_{0}\right\}  $ is a tubular neighborhood of $\Sigma$ and for $\xi,\zeta
\in\Sigma$ with $d\left(  \xi,\zeta\right)  <2\delta_{0}$, we have $d_{\Sigma
}\left(  \xi,\zeta\right)  \leq2d\left(  \xi,\zeta\right)  $. For $x\in V$,
let $\pi\left(  x\right)  \in\Sigma$ be the unique nearest point on $\Sigma$
to $x$. For $\delta>0$, we write $M_{\delta}=\left\{  x\in M:t\left(
x\right)  \leq\delta\right\}  $. For $x\in\overline{M}$, $\delta>0$, we use
$B_{\delta}\left(  x\right)  $ to denote the ball with center at $x$, radius
$\delta$ in $\left(  \overline{M},g\right)  $.

The Green's function of the Laplace operator satisfies%
\[
\left\{
\begin{array}
[c]{l}%
-\Delta_{x}G\left(  x,y\right)  =\delta_{y}\text{ on }M,\\
G\left(  x,y\right)  =0\text{ for }x\in\Sigma.
\end{array}
\right.
\]
Note that $G\left(  x_{1},x_{2}\right)  =G\left(  x_{2},x_{1}\right)  $ for
$x_{1},x_{2}\in M$.

\begin{itemize}
\item The solution of $\left\{
\begin{array}
[c]{l}%
-\Delta u=h\text{ on }M,\\
\left.  u\right\vert _{\Sigma}=0
\end{array}
\right.  $ is given by%
\[
u\left(  x\right)  =\int_{M}G\left(  x,y\right)  h\left(  y\right)
d\mu\left(  y\right)  .
\]

\item The solution of $\left\{
\begin{array}
[c]{l}%
-\Delta u=0\text{ on }M,\\
\left.  u\right\vert _{\Sigma}=f
\end{array}
\right.  $ is given by%
\[
u\left(  x\right)  =-\int_{\Sigma}\left.  \frac{\partial G\left(  x,y\right)
}{\partial_{y}\nu}\right\vert _{y=\xi}f\left(  \xi\right)  dS\left(
\xi\right)  .
\]
Here $\nu$ is the unit outer normal direction on $\Sigma$. In particular the
Poisson kernel is given by%
\[
P\left(  x,\xi\right)  =-\left.  \frac{\partial G\left(  x,y\right)
}{\partial_{y}\nu}\right\vert _{y=\xi}.
\]

\item If $\left\{
\begin{array}
[c]{l}%
-\Delta u=h\text{ on }M\\
\left.  u\right\vert _{\Sigma}=0
\end{array}
\right.  $, then $\frac{\partial u}{\partial\nu}\left(  \xi\right)  =-\int
_{M}P\left(  x,\xi\right)  h\left(  x\right)  d\mu\left(  x\right)  $. In the
future we will denote
\[
\left(  Th\right)  \left(  \xi\right)  =\int_{M}P\left(  x,\xi\right)
h\left(  x\right)  d\mu\left(  x\right)  .
\]
Hence $\frac{\partial u}{\partial\nu}=-Th$.
\end{itemize}

For $f$ defined on $\Sigma$, we write%
\[
\left(  Pf\right)  \left(  x\right)  =\int_{\Sigma}P\left(  x,\xi\right)
f\left(  \xi\right)  dS\left(  \xi\right)  .
\]
$Pf$ is the harmonic extension of $f$.

\begin{lemma}
\label{lem2.1}For $0\leq\delta<\delta_{0}$, denote $\Sigma_{\delta}=\left\{
x\in M:d\left(  x,\Sigma\right)  =\delta\right\}  $. If $u\in C^{\infty
}\left(  M\right)  $ is a nonnegative harmonic function, then%
\[
\int_{\Sigma_{\delta}}udS\leq c\left(  M,g\right)  \int_{\Sigma}udS.
\]

\end{lemma}

\begin{proof}
Denote $\nu$ as the unit outer normal direction. Since $\delta_{0}$ is small,
for $0\leq\delta<\delta_{0}$, the map $\psi_{\delta}:\Sigma\rightarrow
\Sigma_{\delta}$ given by $\psi_{\delta}\left(  \xi\right)  =\exp_{\xi}\left(
-\delta\nu\left(  \xi\right)  \right)  $ is a diffeomorphism and $\int
_{\Sigma_{\delta}}udS=\int_{\Sigma}u\circ\psi_{\delta}\cdot J_{\psi_{\delta}%
}dS$. Hence%
\begin{align*}
\frac{d}{d\delta}\int_{\Sigma_{\delta}}udS  & =\int_{\Sigma_{\delta}}%
\frac{\partial u}{\partial t}dS+\int_{\Sigma}u\circ\psi_{\delta}\cdot
\frac{dJ_{\psi_{\delta}}}{d\delta}dS\\
& \leq c\left(  M,g\right)  \int_{\Sigma}u\circ\psi_{\delta}\cdot
J_{\psi_{\delta}}dS\\
& =c\left(  M,g\right)  \int_{\Sigma_{\delta}}udS.
\end{align*}
Here we have used the equation $\int_{\Sigma_{\delta}}\frac{\partial
u}{\partial t}dS=0$ which follows from the divergence theorem and the fact $u$
is harmonic. It follows that $\int_{\Sigma_{\delta}}udS\leq c\left(
M,g\right)  \int_{\Sigma}udS$.
\end{proof}

To avoid confusion we emphasize that the constants $c\left(  M,g\right)  $'s
are different in different formulas. This convention applies throughout the
article. We will need the following classical estimate for Poisson kernels.

\begin{lemma}
\label{lem2.2}The Poisson kernel $P\left(  x,\xi\right)  $ satisfies%
\[
0\leq P\left(  x,\xi\right)  \leq c\left(  M,g\right)  \frac{t\left(
x\right)  }{\left[  t\left(  x\right)  ^{2}+d_{\Sigma}\left(  \pi\left(
x\right)  ,\xi\right)  ^{2}\right]  ^{\frac{n}{2}}}%
\]
for $x\in M_{\delta_{0}}$ and $\xi\in\Sigma$.
\end{lemma}

\begin{proof}
It follows from Lemma \ref{lem2.1} and an approximation procedure that for
$0<\delta\leq\delta_{0}$,%
\[
\int_{M_{\delta}}P\left(  x,\xi\right)  d\mu\left(  x\right)  \leq c\left(
M,g\right)  \delta.
\]
Since $P\left(  x,\xi\right)  $ is nonnegative, harmonic in $x$ and $P\left(
x,\xi\right)  =0$ for $x\in\Sigma\backslash\left\{  \xi\right\}  $, it follows
from the elliptic estimates of harmonic function that we only need to consider
the case $t\left(  x\right)  +d_{\Sigma}\left(  \pi\left(  x\right)
,\xi\right)  $ is small. Let $t\left(  x\right)  +d_{\Sigma}\left(  \pi\left(
x\right)  ,\xi\right)  =\delta$. If $t\left(  x\right)  \geq\frac{\delta}{7}$,
by mean value inequality%
\begin{align*}
P\left(  x,\xi\right)   & \leq\frac{c\left(  M,g\right)  }{\delta^{n}}%
\int_{B_{\frac{\delta}{7}}\left(  x\right)  }P\left(  y,\xi\right)
d\mu\left(  y\right)  \leq\frac{c\left(  M,g\right)  }{\delta^{n-1}}\\
& \leq c\left(  M,g\right)  \frac{t\left(  x\right)  }{\left[  t\left(
x\right)  ^{2}+d_{\Sigma}\left(  \pi\left(  x\right)  ,\xi\right)
^{2}\right]  ^{\frac{n}{2}}}.
\end{align*}
Assume $t\left(  x\right)  <\frac{\delta}{7}$, then $d\left(  \pi\left(
x\right)  ,\xi\right)  >\frac{3\delta}{7}$. By the gradient estimate of
harmonic functions we know%
\[
\left\vert \nabla P\left(  \cdot,\xi\right)  \right\vert _{L^{\infty}\left(
B_{2\delta/7}\left(  \pi\left(  x\right)  \right)  \cap M\right)  }\leq
\frac{c\left(  M,g\right)  }{\delta^{n+1}}\int_{B_{3\delta/7}\left(
\pi\left(  x\right)  \right)  \cap M}P\left(  y,\xi\right)  d\mu\left(
y\right)  \leq\frac{c\left(  M,g\right)  }{\delta^{n}},
\]
hence $P\left(  x,\xi\right)  \leq c\left(  M,g\right)  \frac{t\left(
x\right)  }{\delta^{n}}$. The lemma follows.
\end{proof}

As an application of Lemma \ref{lem2.2} we may derive the following inequality
for harmonic extensions. Recall if $X$ is a measure space, $p>0$ and $u$ is a
measurable function on $X$, then%
\[
\left\vert u\right\vert _{L_{W}^{p}\left(  X\right)  }=\sup_{t>0}t\left\vert
\left\vert u\right\vert >t\right\vert ^{\frac{1}{p}}.
\]
Here $\left\vert \left\vert u\right\vert >t\right\vert $ is the measure of the
set $\left\{  \left\vert u\right\vert >t\right\}  $.

\begin{proposition}
\label{prop2.1}The harmonic extension operator $P$ satisfies%
\[
\left\vert Pf\right\vert _{L_{W}^{\frac{n}{n-1}}\left(  M\right)  }\leq
c\left(  M,g\right)  \left\vert f\right\vert _{L^{1}\left(  \Sigma\right)  }%
\]
and%
\[
\left\vert Pf\right\vert _{L^{\frac{np}{n-1}}\left(  M\right)  }\leq c\left(
M,g,p\right)  \left\vert f\right\vert _{L^{p}\left(  \Sigma\right)  }%
\]
for $1<p\leq\infty$.
\end{proposition}

\begin{proof}
We only need to prove the weak type estimate. The strong estimate follows from
Marcinkiewicz interpolation theorem (\cite[p197]{SW}) and the basic fact
$\left\vert Pf\right\vert _{L^{\infty}\left(  M\right)  }\leq\left\vert
f\right\vert _{L^{\infty}\left(  \Sigma\right)  }$. To prove the weak type
estimate we may assume $f\geq0$ and $\left\vert f\right\vert _{L^{1}\left(
\Sigma\right)  }=1$. It follows from Lemma \ref{lem2.2} that%
\[
0\leq\left(  Pf\right)  \left(  x\right)  \leq\frac{c\left(  M,g\right)
}{t\left(  x\right)  ^{n-1}}.
\]
For $\delta_{0}=\delta_{0}\left(  M,g\right)  >0$ small, it follows from Lemma
\ref{lem2.1} that
\[
\int_{M_{\delta}}P\left(  x,\xi\right)  d\mu\left(  x\right)  \leq c\left(
M,g\right)  \delta\text{ for }\xi\in\Sigma\text{ and }0<\delta<\delta_{0}.
\]
Hence for $\delta\in\left(  0,\delta_{0}\right)  $,%
\begin{align*}
& \int_{M_{\delta}}\left(  Pf\right)  \left(  x\right)  d\mu\left(  x\right)
\\
& =\int_{\Sigma}dS\left(  \xi\right)  \left[  f\left(  \xi\right)
\int_{M_{\delta}}P\left(  x,\xi\right)  d\mu\left(  x\right)  \right] \\
& \leq c\left(  M,g\right)  \delta.
\end{align*}
For $\lambda\geq c\left(  M,g\right)  $, we have%
\begin{align*}
& \left\vert Pf>\lambda\right\vert \\
& =\left\vert \left\{  x\in M:t\left(  x\right)  <c\left(  M,g\right)
\lambda^{-\frac{1}{n-1}},\left(  Pf\right)  \left(  x\right)  >\lambda
\right\}  \right\vert \\
& \leq\frac{1}{\lambda}\int_{M_{c\left(  M,g\right)  \lambda^{-\frac{1}{n-1}}%
}}\left(  Pf\right)  \left(  x\right)  d\mu\left(  x\right) \\
& \leq c\left(  M,g\right)  \lambda^{-\frac{n}{n-1}}.
\end{align*}
The proposition follows.
\end{proof}

For $1<p<\infty$, if we write%
\begin{equation}
c_{M,g,p}=\sup\left\{  \left\vert Pf\right\vert _{L^{\frac{np}{n-1}}\left(
M\right)  }:f\in L^{p}\left(  \Sigma\right)  ,\left\vert f\right\vert
_{L^{p}\left(  \Sigma\right)  }=1\right\}  ,\label{eq2.3}%
\end{equation}
then $c_{M,g,p}<\infty$. In view of (\ref{eq2.2}), when the background metric
$g$ has zero scalar curvature,%
\begin{equation}
\Theta_{M,g}=c_{M,g,\frac{2\left(  n-1\right)  }{n-2}}^{\frac{2}{n-2}}%
<\infty.\label{eq2.4}%
\end{equation}

In the future we will also need the following compactness property.

\begin{corollary}
\label{cor2.1}For $1\leq p<\infty$, $1\leq q<\frac{np}{n-1}$, the operator
$P:L^{p}\left(  \Sigma\right)  \rightarrow L^{q}\left(  M\right)  $ is compact.
\end{corollary}

\begin{proof}
First assume $1<p<\infty$. If $f_{i}\in L^{p}\left(  \Sigma\right)  $ such
that $\left\vert f_{i}\right\vert _{L^{p}\left(  \Sigma\right)  }\leq1$, it
follows from Lemma \ref{lem2.2} that%
\[
\left\vert \left(  Pf_{i}\right)  \left(  x\right)  \right\vert \leq
\frac{c\left(  M,g\right)  }{t\left(  x\right)  ^{n-1}}\text{ for }x\in
M\backslash\Sigma.
\]
Using elliptic estimates of harmonic functions we know after passing to a
subsequence we may find a $u\in C^{\infty}\left(  M\backslash\Sigma\right)  $
such that $Pf_{i}\rightarrow u$ in $C_{loc}^{\infty}\left(  M\backslash
\Sigma\right)  $. For $\delta>0$ small, we have%
\begin{align*}
& \left\vert Pf_{i}-Pf_{j}\right\vert _{L^{q}\left(  M\right)  }\\
& \leq\left\vert Pf_{i}-Pf_{j}\right\vert _{L^{q}\left(  M\backslash
M_{\delta}\right)  }+\left\vert Pf_{i}-Pf_{j}\right\vert _{L^{q}\left(
M_{\delta}\right)  }\\
& \leq\left\vert Pf_{i}-Pf_{j}\right\vert _{L^{q}\left(  M\backslash
M_{\delta}\right)  }+\left\vert Pf_{i}-Pf_{j}\right\vert _{L^{\frac{np}{n-1}%
}\left(  M_{\delta}\right)  }\left\vert M_{\delta}\right\vert ^{\frac{1}%
{q}-\frac{n-1}{np}}\\
& \leq\left\vert Pf_{i}-Pf_{j}\right\vert _{L^{q}\left(  M\backslash
M_{\delta}\right)  }+c\left(  M,g,p\right)  \left\vert M_{\delta}\right\vert
^{\frac{1}{q}-\frac{n-1}{np}}.
\end{align*}
Hence%
\[
\lim\sup_{i,j\rightarrow\infty}\left\vert Pf_{i}-Pf_{j}\right\vert
_{L^{q}\left(  M\right)  }\leq c\left(  M,g,p\right)  \left\vert M_{\delta
}\right\vert ^{\frac{1}{q}-\frac{n-1}{np}}.
\]
Letting $\delta\rightarrow0^{+}$, we see $Pf_{i}$ is a Cauchy sequence in
$L^{q}\left(  M\right)  $. In another word, $P:L^{p}\left(  \Sigma\right)
\rightarrow L^{q}\left(  M\right)  $ is compact.

When $p=1$, the argument is similar. We only need to observe that for any
$1\leq q<\widetilde{q}<\frac{n}{n-1}$, $P:L^{1}\left(  \Sigma\right)
\rightarrow L^{\widetilde{q}}\left(  M\right)  $ is bounded.
\end{proof}

Let $h$ be a function on $M$, recall $\left(  Th\right)  \left(  \xi\right)
=\int_{M}P\left(  x,\xi\right)  h\left(  x\right)  d\mu\left(  x\right)  $. We
have the following dual statement to Proposition \ref{prop2.1}.

\begin{proposition}
\label{prop2.2}For $1\leq p<n$ and $h\in L^{p}\left(  M\right)  $,%
\[
\left\vert Th\right\vert _{L^{\frac{\left(  n-1\right)  p}{n-p}}\left(
\Sigma\right)  }\leq c\left(  M,g,p\right)  \left\vert h\right\vert
_{L^{p}\left(  M\right)  }.
\]

\end{proposition}

\begin{proof}
We may prove the inequality by a duality argument. Indeed for any nonnegative
functions $h$ on $M$ and $f$ on $\Sigma$, we have%
\begin{align*}
0  & \leq\int_{\Sigma}\left(  Th\right)  \left(  \xi\right)  f\left(
\xi\right)  dS\left(  \xi\right)  =\int_{\Sigma}dS\left(  \xi\right)  \int
_{M}P\left(  x,\xi\right)  h\left(  x\right)  f\left(  \xi\right)  d\mu\left(
x\right) \\
& =\int_{M}\left(  Pf\right)  \left(  x\right)  h\left(  x\right)  d\mu\left(
x\right)  \leq\left\vert Pf\right\vert _{L^{\frac{p}{p-1}}\left(  M\right)
}\left\vert h\right\vert _{L^{p}\left(  M\right)  }\\
& \leq c\left(  M,g,p\right)  \left\vert h\right\vert _{L^{p}\left(  M\right)
}\left\vert f\right\vert _{L^{\frac{\left(  n-1\right)  p}{n\left(
p-1\right)  }}\left(  \Sigma\right)  },
\end{align*}
the proposition follows. One may also prove the inequality directly. Indeed it
follows from Lemma \ref{lem2.2} that $\left\vert P\left(  \cdot,\xi\right)
\right\vert _{L^{\frac{n}{n-1},\infty}\left(  M\right)  }\leq c\left(
M,g\right)  <\infty$ for $\xi\in\Sigma$. Hence $T:L^{n,1}\left(  M\right)
\rightarrow L^{\infty}\left(  \Sigma\right)  $ is a bounded linear map. On the
other hand for $h\in L^{1}\left(  M\right)  $,%
\[
\int_{\Sigma}\left\vert \left(  Th\right)  \left(  \xi\right)  \right\vert
dS\left(  \xi\right)  \leq\int_{\Sigma}dS\left(  \xi\right)  \int_{M}P\left(
x,\xi\right)  \left\vert h\left(  x\right)  \right\vert d\mu\left(  x\right)
=\int_{M}\left\vert h\left(  x\right)  \right\vert d\mu\left(  x\right)  .
\]
Hence $T:L^{1}\left(  M\right)  \rightarrow L^{1}\left(  \Sigma\right)  $ is
also bounded. The proposition follows from the Marcinkiewicz interpolation
theorem. Finally we point out for $1<p<n$, we may solve $\left\{
\begin{array}
[c]{l}%
-\Delta u=h\text{ on }M\\
\left.  u\right\vert _{\Sigma}=0
\end{array}
\right.  $ and $\left(  Th\right)  \left(  \xi\right)  =-\frac{\partial
u}{\partial\nu}\left(  \xi\right)  $. By the $L^{p}$ theory we know
$\left\vert u\right\vert _{W^{2,p}\left(  M\right)  }\leq c\left(
M,g,p\right)  \left\vert h\right\vert _{L^{p}\left(  M\right)  }$. It follows
from boundary trace imbedding theorem (\cite[p164]{A}) that%
\[
\left\vert Th\right\vert _{L^{\frac{\left(  n-1\right)  p}{n-p}}\left(
\Sigma\right)  }=\left\vert \frac{\partial u}{\partial\nu}\right\vert
_{L^{\frac{\left(  n-1\right)  p}{n-p}}\left(  \Sigma\right)  }\leq c\left(
M,g,p\right)  \left\vert u\right\vert _{W^{2,p}\left(  M\right)  }\leq
c\left(  M,g,p\right)  \left\vert h\right\vert _{L^{p}\left(  M\right)  }.
\]

\end{proof}

\subsection{Miscellaneous\label{sec2.2}}

Later on we will need the following Hausdorff-Young type inequality to
estimate some nonmajor terms.

\begin{lemma}
\label{lem2.3}Let $X$ and $Y$ be measure spaces, $1\leq p,q_{0},q_{1}%
,r\leq\infty$, $p\leq r$, $q_{0}\leq r$ and%
\[
\frac{1}{p}+\frac{1}{q_{1}}=\frac{q_{0}}{q_{1}r}+1.
\]
Assume $K$ is defined on $X\times Y$ such that%
\[
\left(  \int_{X}\left\vert K\left(  x,y\right)  \right\vert ^{q_{0}}dx\right)
^{\frac{1}{q_{0}}}\leq A,\quad\left(  \int_{Y}\left\vert K\left(  x,y\right)
\right\vert ^{q_{1}}dy\right)  ^{\frac{1}{q_{1}}}\leq A.
\]
For a function $f$ defined on $Y$, we let $\left(  Kf\right)  \left(
x\right)  =\int_{Y}K\left(  x,y\right)  f\left(  y\right)  dy$, then%
\[
\left\vert Kf\right\vert _{L^{r}\left(  X\right)  }\leq A\left\vert
f\right\vert _{L^{p}\left(  Y\right)  }.
\]

\end{lemma}

\begin{proof}
Without losing of generality we may assume $K\geq0$ and $f\geq0$, then%
\begin{align*}
& \left(  Kf\right)  \left(  x\right) \\
& =\int_{Y}K\left(  x,y\right)  ^{\frac{q_{0}}{r}}f\left(  y\right)
^{\frac{p}{r}}K\left(  x,y\right)  ^{\frac{r-q_{0}}{r}}f\left(  y\right)
^{\frac{r-p}{r}}dy\\
& \leq\left(  \int_{Y}K\left(  x,y\right)  ^{q_{0}}f\left(  y\right)
^{p}dy\right)  ^{\frac{1}{r}}\left(  \int_{Y}K\left(  x,y\right)  ^{q_{1}%
}dy\right)  ^{\frac{r-q_{0}}{q_{1}r}}\left(  \int_{Y}f\left(  y\right)
^{p}dy\right)  ^{\frac{r-p}{pr}}\\
& \leq A^{\frac{r-q_{0}}{r}}\left\vert f\right\vert _{L^{p}\left(  Y\right)
}^{\frac{r-p}{r}}\left(  \int_{Y}K\left(  x,y\right)  ^{q_{0}}f\left(
y\right)  ^{p}dy\right)  ^{\frac{1}{r}}.
\end{align*}
Here we have used the Holder's inequality and the fact $\frac{1}{r}+\frac
{1}{\frac{q_{1}r}{r-q_{0}}}+\frac{1}{\frac{pr}{r-p}}=1$. Hence%
\[
\left(  Kf\right)  \left(  x\right)  ^{r}\leq A^{r-q_{0}}\left\vert
f\right\vert _{L^{p}\left(  Y\right)  }^{r-p}\int_{Y}K\left(  x,y\right)
^{q_{0}}f\left(  y\right)  ^{p}dy.
\]
Integrating both sides, we get the needed inequality.
\end{proof}

\section{Sharp inequalities on the unit ball\label{sec3}}

The aim of this section is to show $\Theta_{\overline{B}_{1},g_{\mathbb{R}%
^{n}}}=c_{\overline{B}_{1},g_{\mathbb{R}^{n}},\frac{2\left(  n-1\right)
}{n-2}}^{\frac{2}{n-2}}=I\left(  \overline{B}_{1},g_{\mathbb{R}^{n}}\right)
=n^{-\frac{1}{n-1}}\omega_{n}^{-\frac{1}{n\left(  n-1\right)  }}$ (see
(\ref{eq2.3}), (\ref{eq2.4})).

\begin{theorem}
\label{thm3.1}Assume $n\geq3$, then for every $f\in L^{\frac{2\left(
n-1\right)  }{n-2}}\left(  \partial B_{1}^{n}\right)  $,%
\[
\left\vert Pf\right\vert _{L^{\frac{2n}{n-2}}\left(  B_{1}\right)  }\leq
n^{-\frac{n-2}{2\left(  n-1\right)  }}\omega_{n}^{-\frac{n-2}{2n\left(
n-1\right)  }}\left\vert f\right\vert _{L^{\frac{2\left(  n-1\right)  }{n-2}%
}\left(  \partial B_{1}\right)  }.
\]
Here $Pf$ is the harmonic extension of $f$, $\omega_{n}$ is the volume of the
unit ball in $\mathbb{R}^{n}$. Equality holds if and only if $f\left(
\xi\right)  =c\left(  1+\lambda\xi\cdot\zeta\right)  ^{-\frac{n-2}{2}}$ for
some constant $c$, $\zeta\in\partial B_{1}$ and $0\leq\lambda<1$.
\end{theorem}

Note that this theorem is a consequence of \cite[Theorem 1.1]{HWY} (see the
discussions before \cite[Theorem 1.1]{HWY}). Below we will present a different
argument which has its own interest. Before discussing the approach, we
describe some corollaries of the theorem. Note that in Proposition
\ref{prop2.1} the strong inequality is not true for $p=1$. Instead we have the following

\begin{corollary}
\label{cor3.1}Assume $n\geq3$, then for $f\in L^{\infty}\left(  \partial
B_{1}^{n}\right)  $,%
\[
\left\vert e^{Pf}\right\vert _{L^{\frac{n}{n-1}}\left(  B_{1}^{n}\right)
}\leq n^{-1}\omega_{n}^{-\frac{1}{n}}\left\vert e^{f}\right\vert
_{L^{1}\left(  \partial B_{1}\right)  }.
\]
Moreover equality holds if and only if $f$ is constant.
\end{corollary}

\begin{proof}
If $u$ is a harmonic function, then $\Delta e^{u}=e^{u}\left\vert \nabla
u\right\vert ^{2}$. Hence $e^{u}$ is subharmonic and not harmonic except when
$u$ is a constant function. It follows from Theorem \ref{thm3.1} that%
\begin{align*}
\left\vert e^{\frac{n-2}{2\left(  n-1\right)  }Pf}\right\vert _{L^{\frac
{2n}{n-2}}\left(  B_{1}\right)  }  & \leq\left\vert P\left(  e^{\frac
{n-2}{2\left(  n-1\right)  }f}\right)  \right\vert _{L^{\frac{2n}{n-2}}\left(
B_{1}\right)  }\\
& \leq n^{-\frac{n-2}{2\left(  n-1\right)  }}\omega_{n}^{-\frac{n-2}{2n\left(
n-1\right)  }}\left\vert e^{\frac{n-2}{2\left(  n-1\right)  }f}\right\vert
_{L^{\frac{2\left(  n-1\right)  }{n-2}}\left(  \partial B_{1}\right)  }.
\end{align*}
Hence%
\[
\left\vert e^{Pf}\right\vert _{L^{\frac{n}{n-1}}\left(  B_{1}\right)  }\leq
n^{-1}\omega_{n}^{-\frac{1}{n}}\left\vert e^{f}\right\vert _{L^{1}\left(
\partial B_{1}\right)  }.
\]
If equality holds, then $e^{\frac{n-2}{2\left(  n-1\right)  }Pf}=P\left(
e^{\frac{n-2}{2\left(  n-1\right)  }f}\right)  $ and $e^{\frac{n-2}{2\left(
n-1\right)  }Pf}$ must be a harmonic function, hence $Pf$ is equal to constant
and so is $f$.
\end{proof}

\begin{corollary}
\label{cor3.2}Assume $n\geq3$, then for $\frac{2\left(  n-1\right)  }%
{n-2}<p<\infty$, $f\in L^{p}\left(  \partial B_{1}^{n}\right)  $,%
\[
\left\vert Pf\right\vert _{L^{\frac{np}{n-1}}\left(  B_{1}\right)  }\leq
n^{-\frac{1}{p}}\omega_{n}^{-\frac{1}{np}}\left\vert f\right\vert
_{L^{p}\left(  \partial B_{1}\right)  }.
\]
Equality holds if and only if $f$ is constant.
\end{corollary}

\begin{proof}
Denote $r=\frac{p}{\frac{2\left(  n-1\right)  }{n-2}}>1$. If $u$ is a harmonic
function on $B_{1}$, then $\left\vert u\right\vert ^{r}$ is a subharmonic
function and it is not harmonic except when $u$ is a constant function. If
$f\in L^{p}\left(  \partial B_{1}\right)  $, then by Theorem \ref{thm3.1},%
\[
\left\vert \left\vert Pf\right\vert ^{r}\right\vert _{L^{\frac{2n}{n-2}%
}\left(  B_{1}\right)  }\leq\left\vert P\left(  \left\vert f\right\vert
^{r}\right)  \right\vert _{L^{\frac{2n}{n-2}}\left(  B_{1}\right)  }\leq
n^{-\frac{n-2}{2\left(  n-1\right)  }}\omega_{n}^{-\frac{n-2}{2n\left(
n-1\right)  }}\left\vert \left\vert f\right\vert ^{r}\right\vert
_{L^{\frac{2\left(  n-1\right)  }{n-2}}\left(  \partial B_{1}\right)  }.
\]
Hence%
\[
\left\vert Pf\right\vert _{L^{\frac{np}{n-1}}\left(  B_{1}\right)  }\leq
n^{-\frac{1}{p}}\omega_{n}^{-\frac{1}{np}}\left\vert f\right\vert
_{L^{p}\left(  \partial B_{1}\right)  }.
\]
If equality holds then $\left\vert Pf\right\vert ^{r}=P\left(  \left\vert
f\right\vert ^{r}\right)  $. In particular $\left\vert Pf\right\vert ^{r}$ is
a harmonic function and hence $Pf$ is a constant function, so is $f$.
\end{proof}

\begin{remark}
\label{rmk3.1}For $1<p<\frac{2\left(  n-1\right)  }{n-2}$, $1$ is still a
critical point for the functional $\frac{\left\vert Pf\right\vert
_{L^{\frac{np}{n-1}}\left(  B_{1}\right)  }}{\left\vert f\right\vert
_{L^{p}\left(  \partial B_{1}\right)  }}$, but calculation shows for
$f_{\varepsilon}\left(  \xi\right)  =1+\varepsilon\xi_{1}$,%
\[
\frac{\left\vert Pf_{\varepsilon}\right\vert _{L^{\frac{np}{n-1}}\left(
B_{1}\right)  }}{\left\vert f_{\varepsilon}\right\vert _{L^{p}\left(  \partial
B_{1}\right)  }}=n^{-\frac{1}{p}}\omega_{n}^{-\frac{1}{np}}\left[
1+\frac{n-2}{2n\left(  n-1\right)  \left(  n+2\right)  }\left(  \frac{2\left(
n-1\right)  }{n-2}-p\right)  \varepsilon^{2}+O\left(  \varepsilon^{4}\right)
\right]  .
\]
Hence $1$ is not a local maximizer. It remains an interesting question to
calculate%
\[
\sup\left\{  \left\vert Pf\right\vert _{L^{\frac{np}{n-1}}\left(
B_{1}\right)  }:f\in L^{p}\left(  \partial B_{1}\right)  ,\left\vert
f\right\vert _{L^{p}\left(  \partial B_{1}\right)  }=1\right\}
\]
for these $p^{\prime}s$.
\end{remark}

The new approach to Theorem \ref{thm3.1} needs an interesting Kazdan-Warner
type condition. To formulate the condition, we introduce the weighted
isoperimetric ratio.

Assume $n\geq2$, $\left(  M^{n},g\right)  $ is a smooth compact Riemannian
manifold with boundary $\Sigma=\partial M$. Let $K$ be a positive smooth
function on $\Sigma$, then we write the weighted isoperimetric ratio%
\[
I\left(  M,g,K\right)  =\frac{\mu\left(  M\right)  ^{\frac{1}{n}}}{\left(
\int_{\Sigma}KdS\right)  ^{\frac{1}{n-1}}}.
\]
Here $d\mu$ is the measure associated with $g$ and $dS$ is the measure on
$\Sigma$. If $n\geq3$ and $\left(  M^{n},g\right)  $ satisfies $\lambda
_{1}\left(  L_{g}\right)  >0$, for $\widetilde{g}\in\left[  g\right]  $ with
zero scalar curvature, we write $\widetilde{g}=u^{\frac{4}{n-2}}g$, $\left.
u\right\vert _{\Sigma}=f$, then%
\[
I\left(  M,\widetilde{g},K\right)  =\frac{\left(  \int_{M}\left(
P_{L}f\right)  ^{\frac{2n}{n-2}}d\mu\right)  ^{\frac{1}{n}}}{\left(
\int_{\Sigma}Kf^{\frac{2\left(  n-1\right)  }{n-2}}dS\right)  ^{\frac{1}{n-1}%
}}.
\]
The Euler-Lagrange equation of this functional reads as%
\[
\int_{M}P_{L}\left(  x,\xi\right)  \left(  P_{L}f\right)  \left(  x\right)
^{\frac{n+2}{n-2}}d\mu\left(  x\right)  =\operatorname*{const}\cdot K\left(
\xi\right)  f\left(  \xi\right)  ^{\frac{n}{n-2}}.
\]

\begin{lemma}
[Kazdan-Warner type condition]\label{lem3.1}Assume $n\geq3$, $\left(
M^{n},g\right)  $ is a smooth compact Riemannian manifold with boundary and
$\lambda_{1}\left(  L_{g}\right)  >0$, $K$ and $f$ are positive smooth
functions on $\Sigma$ such that%
\[
\int_{M}P_{L}\left(  x,\xi\right)  \left(  P_{L}f\right)  \left(  x\right)
^{\frac{n+2}{n-2}}d\mu\left(  x\right)  =K\left(  \xi\right)  f\left(
\xi\right)  ^{\frac{n}{n-2}}.
\]
Let $X$ be a conformal vector field on $M$ (note $X$ must be tangent to
$\Sigma$), then%
\[
\int_{\Sigma}XK\cdot f^{\frac{2\left(  n-1\right)  }{n-2}}dS=0.
\]

\end{lemma}

\begin{proof}
Denote $u=P_{L}f$. Let $\phi_{t}$ be the smooth $1$-parameter group generated
by $X$, then%
\[
\left.  \frac{d}{dt}\right\vert _{t=0}I\left(  M,\phi_{t}^{\ast}\left(
u^{\frac{4}{n-2}}g\right)  ,K\right)  =0.
\]
On the other hand,%
\begin{align*}
\left.  \frac{d}{dt}\right\vert _{t=0}I\left(  M,\phi_{t}^{\ast}\left(
u^{\frac{4}{n-2}}g\right)  ,K\right)   & =\left.  \frac{d}{dt}\right\vert
_{t=0}I\left(  M,u^{\frac{4}{n-2}}g,K\circ\phi_{-t}\right) \\
& =\frac{I\left(  M,u^{\frac{4}{n-2}}g,K\right)  }{n-1}\frac{\int_{\Sigma
}XK\cdot f^{\frac{2\left(  n-1\right)  }{n-2}}dS}{\int_{\Sigma}Kf^{\frac
{2\left(  n-1\right)  }{n-2}}dS}.
\end{align*}
This implies $\int_{\Sigma}XK\cdot f^{\frac{2\left(  n-1\right)  }{n-2}}dS=0$.
\end{proof}

\begin{corollary}
\label{cor3.3}Assume $n\geq3$, $K$ and $f$ are positive smooth functions on
$\partial B_{1}^{n}$ such that%
\[
\int_{B_{1}}P\left(  x,\xi\right)  \left(  Pf\right)  \left(  x\right)
^{\frac{n+2}{n-2}}dx=K\left(  \xi\right)  f\left(  \xi\right)  ^{\frac{n}%
{n-2}},
\]
then $\int_{\partial B_{1}}\left\langle \nabla K\left(  \xi\right)  ,\nabla
\xi_{i}\right\rangle f\left(  \xi\right)  ^{\frac{2\left(  n-1\right)  }{n-2}%
}dS\left(  \xi\right)  =0$ for $1\leq i\leq n$.
\end{corollary}

This is because $\nabla\xi_{i}$ is the restriction to $\partial B_{1}$ of a
conformal vector field on $\left(  \overline{B}_{1},g_{\mathbb{R}^{n}}\right)
$.

We will also need some rearrangement inequality on $\partial B_{1}$ which was
proven in \cite{BT}. We say a function $f$ on $\partial B_{1}$ is radially
symmetric if $f\left(  \xi\right)  $ is a function of $\xi_{n}$. Let $f$ be a
measurable function on $\partial B_{1}$, then the symmetric rearrangement of
$f$ is a radial decreasing function $f^{\ast}$ which has the same distribution
as $f$. The following rearrangement inequality was proven in \cite[Theorem
2]{BT}. Namely, if $K$ is a nondecreasing bounded function on $\left[
-1,1\right]  $, then for all $f,g\in L^{1}\left(  \partial B_{1}\right)  $,%
\begin{align*}
& \int_{\partial B_{1}\times\partial B_{1}}f\left(  \xi\right)  g\left(
\eta\right)  K\left(  \xi\cdot\eta\right)  dS\left(  \xi\right)  dS\left(
\eta\right) \\
& \leq\int_{\partial B_{1}\times\partial B_{1}}f^{\ast}\left(  \xi\right)
g^{\ast}\left(  \eta\right)  K\left(  \xi\cdot\eta\right)  dS\left(
\xi\right)  dS\left(  \eta\right)  .
\end{align*}
It follows that if $K$ is a bounded nonnegative nondecreasing function on
$\left[  -1,1\right]  $, $f$ is nonnegative function on $\partial B_{1}$ and%
\[
\left(  K\ast f\right)  \left(  \xi\right)  =\int_{\partial B_{1}}K\left(
\xi\cdot\eta\right)  f\left(  \eta\right)  dS\left(  \eta\right)  ,
\]
then for $1\leq p<\infty$, $\left\vert K\ast f\right\vert _{L^{p}\left(
\partial B_{1}\right)  }\leq\left\vert K\ast f^{\ast}\right\vert
_{L^{p}\left(  \partial B_{1}\right)  }$.

Recall the Poisson kernel on $\left(  \overline{B}_{1},g_{\mathbb{R}^{n}%
}\right)  $ is given by%
\[
P\left(  x,\xi\right)  =\frac{1-\left\vert x\right\vert ^{2}}{n\omega
_{n}\left\vert x-\xi\right\vert ^{n}}.
\]
For $0<r<1$, $\xi,\zeta\in\partial B_{1}$,%
\[
P\left(  r\zeta,\xi\right)  =\frac{1-r^{2}}{n\omega_{n}\left(  r^{2}%
+1-2r\zeta\cdot\xi\right)  ^{\frac{n}{2}}}=K_{r}\left(  \zeta\cdot\xi\right)
.
\]
Hence for $1\leq p<\infty$ and $f\geq0$,%
\begin{align*}
\left\vert Pf\right\vert _{L^{p}\left(  B_{1}\right)  }^{p}  & =\int_{0}%
^{1}\left\vert K_{r}\ast f\right\vert _{L^{p}\left(  \partial B_{1}\right)
}^{p}r^{n-1}dr\\
& \leq\int_{0}^{1}\left\vert K_{r}\ast f^{\ast}\right\vert _{L^{p}\left(
\partial B_{1}\right)  }^{p}r^{n-1}dr=\left\vert Pf^{\ast}\right\vert
_{L^{p}\left(  B_{1}\right)  }^{p}.
\end{align*}
It follows that $\left\vert Pf\right\vert _{L^{p}\left(  B_{1}\right)  }%
\leq\left\vert Pf^{\ast}\right\vert _{L^{p}\left(  B_{1}\right)  }$.

\begin{proof}
[Proof of Theorem \ref{thm3.1}]For $p>\frac{2\left(  n-1\right)  }{n-2}$, we
consider the variational problem%
\begin{equation}
\sup\left\{  \left\vert Pf\right\vert _{L^{\frac{2n}{n-2}}\left(
B_{1}\right)  }:f\in L^{p}\left(  \partial B_{1}\right)  ,\left\vert
f\right\vert _{L^{p}\left(  \partial B_{1}\right)  }=1\right\}  .\label{eq3.1}%
\end{equation}
By Corollary \ref{cor2.1} the operator $P:L^{p}\left(  \partial B_{1}\right)
\rightarrow L^{\frac{2n}{n-2}}\left(  B_{1}\right)  $ is compact, hence the
supreme is achieved at some $f_{p}\geq0$. Replacing $f_{p}$ by $f_{p}^{\ast}%
$\ we may assume $f_{p}$ is radial symmetric and decreasing. After scaling
$f_{p}$ satisfies%
\[
f_{p}\left(  \xi\right)  ^{p-1}=\int_{B_{1}}P\left(  x,\xi\right)  \left(
Pf_{p}\right)  \left(  x\right)  ^{\frac{n+2}{n-2}}dx.
\]
Standard bootstrap using Proposition \ref{prop2.1} and Proposition
\ref{prop2.2} shows $f_{p}\in C^{\infty}\left(  \partial B_{1}\right)  $ and
$f_{p}>0$. Rewrite the equation as%
\[
\int_{B_{1}}P\left(  x,\xi\right)  \left(  Pf_{p}\right)  \left(  x\right)
^{\frac{n+2}{n-2}}dx=f_{p}\left(  \xi\right)  ^{\frac{n}{n-2}}f_{p}\left(
\xi\right)  ^{p-\frac{2\left(  n-1\right)  }{n-2}}.
\]
It follows from Corollary \ref{cor3.3} that%
\[
\int_{\partial B_{1}}\left\langle \nabla f_{p}\left(  \xi\right)
^{p-\frac{2\left(  n-1\right)  }{n-2}},\nabla\xi_{n}\right\rangle f_{p}\left(
\xi\right)  ^{\frac{2\left(  n-1\right)  }{n-2}}dS\left(  \xi\right)  =0.
\]
We may write $g_{p}\left(  r\right)  =f_{p}\left(  0,\cdots,0,\sin r,\cos
r\right)  $ for $0\leq r\leq\pi$. Then the equality becomes $\int_{0}^{\pi
}g_{p}^{\prime}\left(  r\right)  g_{p}\left(  r\right)  ^{p-1}\sin^{n-1}%
rdr=0$. Since $g_{p}^{\prime}\leq0$ and $g_{p}>0$, we get $g_{p}^{\prime}=0$
and hence $f_{p}\equiv\operatorname*{const}$. This implies%
\[
\left\vert Pf\right\vert _{L^{\frac{2n}{n-2}}\left(  B_{1}\right)  }\leq
\frac{\omega_{n}^{\frac{n-2}{2n}}}{\left(  n\omega_{n}\right)  ^{\frac{1}{p}}%
}\left\vert f\right\vert _{L^{p}\left(  \partial B_{1}\right)  }.
\]
Let $p\rightarrow\frac{2\left(  n-1\right)  }{n-2}$, we get the needed
inequality. At last we may apply \cite[Theorem 1.2]{HWY} to identify all the
functions which achieves the equality.
\end{proof}

\section{Regularity of solutions to some nonlinear integral
equations\label{sec4}}

Assume $1<p<\infty$. If $f\in L^{p}\left(  \Sigma\right)  $ is a maximizer for
the variational problem%
\[
c_{M,g,p}=\sup\left\{  \left\vert Pf\right\vert _{L^{\frac{np}{n-1}}\left(
M\right)  }:f\in L^{p}\left(  \Sigma\right)  ,\left\vert f\right\vert
_{L^{p}\left(  \Sigma\right)  }=1\right\}  ,
\]
then we may assume $f\geq0$, moreover after suitable scaling it satisfies the
nonlinear integral equation%
\[
f\left(  \xi\right)  ^{p-1}=\int_{M}P\left(  x,\xi\right)  \left(  Pf\right)
\left(  x\right)  ^{\frac{np}{n-1}-1}d\mu\left(  x\right)  .
\]
This section is aiming at proving all these solutions are in fact smooth.

\begin{proposition}
\label{prop4.1}Assume $n\geq2$, $\left(  M^{n},g\right)  $ is a smooth compact
Riemannian manifold with boundary $\Sigma=\partial M$. If $1<p<\infty$, $f\in
L^{p}\left(  \Sigma\right)  $ is nonnegative, not identically zero and it
satisfies%
\[
f\left(  \xi\right)  ^{p-1}=\int_{M}P\left(  x,\xi\right)  \left(  Pf\right)
\left(  x\right)  ^{\frac{np}{n-1}-1}d\mu\left(  x\right)  ,
\]
then $f\in C^{\infty}\left(  \Sigma\right)  $.
\end{proposition}

\begin{proof}
Let $p_{0}=\frac{1}{p-1}$, $f_{0}\left(  \xi\right)  =f\left(  \xi\right)
^{p-1}$, $u_{0}\left(  x\right)  =\left(  Pf\right)  \left(  x\right)  $, then
$0<p_{0}<\infty$, $f_{0}\in L^{p_{0}+1}\left(  \Sigma\right)  $, $u_{0}\in
L^{\frac{n\left(  p_{0}+1\right)  }{\left(  n-1\right)  p_{0}}}\left(
M\right)  $ and%
\[
u_{0}\left(  x\right)  =\int_{\Sigma}P\left(  x,\xi\right)  f_{0}\left(
\xi\right)  ^{p_{0}}dS\left(  \xi\right)  ,\quad f_{0}\left(  \xi\right)
=\int_{M}P\left(  x,\xi\right)  u_{0}\left(  x\right)  ^{\frac{p_{0}%
+n}{\left(  n-1\right)  p_{0}}}d\mu\left(  x\right)  .
\]
Let $\left(  \overline{M},g\right)  $ be the same as in Section \ref{sec2.1}.
Given $\xi_{0}\in\Sigma$, by choosing a local coordinate $\phi:U\left(
\xi_{0}\right)  \rightarrow\left\{  x\in\mathbb{R}^{n}:\left\vert x\right\vert
<2\right\}  $ with $\phi\left(  \xi_{0}\right)  =0$ and $\phi\left(  U\left(
\xi_{0}\right)  \cap M\right)  =\left\{  x\in\mathbb{R}^{n}:\left\vert
x\right\vert <2,x_{n}\geq0\right\}  $, we may identify $U\left(  \xi
_{0}\right)  $ with $\left\{  x\in\mathbb{R}^{n}:\left\vert x\right\vert
<2\right\}  $. For $0<R<1$, we write
\begin{align*}
B_{R}^{+}  & =\left\{  x\in\mathbb{R}^{n}:\left\vert x\right\vert
<R,x_{n}>0\right\}  ,\\
B_{R}  & =B_{R}^{n-1}=\left\{  \xi\in\mathbb{R}^{n-1}:\left\vert
\xi\right\vert <R\right\}
\end{align*}
and%
\begin{align*}
u_{R}\left(  x\right)   & =\int_{\Sigma\backslash B_{R}}P\left(  x,\xi\right)
f_{0}\left(  \xi\right)  ^{p_{0}}dS\left(  \xi\right)  ,\\
f_{R}\left(  \xi\right)   & =\int_{M\backslash B_{R}^{+}}P\left(
x,\xi\right)  u_{0}\left(  x\right)  ^{\frac{p_{0}+n}{\left(  n-1\right)
p_{0}}}d\mu\left(  x\right)  .
\end{align*}
Then $u_{R}\in C^{\infty}\left(  \left\{  x\in\mathbb{R}^{n}:\left\vert
x\right\vert <R,x_{n}\geq0\right\}  \right)  $, $f_{R}\in C^{\infty}\left(
B_{R}\right)  $. To prove the regularity of $f$, we discuss two cases.

\begin{case}
\label{case4.1}$0<p_{0}\leq\frac{n}{n-1}$.
\end{case}

In this case, we have $\frac{p_{0}+n}{\left(  n-1\right)  p_{0}}>1$. Fix a
number $r$ such that%
\[
1\leq r<\frac{p_{0}+n}{\left(  n-1\right)  p_{0}}\text{ and }r>\frac{1}{p_{0}%
},
\]
then%
\[
f_{0}\left(  \xi\right)  ^{1/r}\leq\left[  \int_{B_{R}^{+}}P\left(
x,\xi\right)  u_{0}\left(  x\right)  ^{\frac{p_{0}+n}{\left(  n-1\right)
p_{0}}}d\mu\left(  x\right)  \right]  ^{1/r}+f_{R}\left(  \xi\right)  ^{1/r}.
\]
Hence using Lemma \ref{lem2.2} we have%
\begin{align*}
& u_{0}\left(  x\right) \\
& =\int_{B_{R}}P\left(  x,\xi\right)  f_{0}\left(  \xi\right)  ^{p_{0}-r^{-1}%
}f_{0}\left(  \xi\right)  ^{1/r}dS\left(  \xi\right)  +u_{R}\left(  x\right)
\\
& \leq\int_{B_{R}}P\left(  x,\xi\right)  f_{0}\left(  \xi\right)
^{p_{0}-r^{-1}}\left[  \int_{B_{R}^{+}}P\left(  y,\xi\right)  u_{0}\left(
y\right)  ^{\frac{p_{0}+n}{\left(  n-1\right)  p_{0}}-r}u_{0}\left(  y\right)
^{r}d\mu\left(  y\right)  \right]  ^{1/r}dS\left(  \xi\right) \\
& +v_{R}\left(  x\right) \\
& \leq c\left(  M,g,p,r\right)  \int_{B_{R}}\frac{x_{n}}{\left(  \left\vert
x^{\prime}-\xi\right\vert ^{2}+x_{n}^{2}\right)  ^{n/2}}f_{0}\left(
\xi\right)  ^{p_{0}-r^{-1}}\cdot\\
& \left[  \int_{B_{R}^{+}}\frac{y_{n}}{\left(  \left\vert y^{\prime}%
-\xi\right\vert ^{2}+y_{n}^{2}\right)  ^{n/2}}u_{0}\left(  y\right)
^{\frac{p_{0}+n}{\left(  n-1\right)  p_{0}}-r}u_{0}\left(  y\right)
^{r}dy\right]  ^{1/r}d\xi+v_{R}\left(  x\right)
\end{align*}
here $dx$ and $d\xi$ means the standard Lebesgue measure and%
\[
v_{R}\left(  x\right)  =\int_{B_{R}}P\left(  x,\xi\right)  f_{0}\left(
\xi\right)  ^{p_{0}-r^{-1}}f_{R}\left(  \xi\right)  ^{1/r}dS\left(
\xi\right)  +u_{R}\left(  x\right)  .
\]
We have $v_{R}\in L^{\frac{n\left(  p_{0}+1\right)  }{\left(  n-1\right)
p_{0}}}\left(  B_{R}^{+}\right)  \cap L_{loc}^{\frac{n\left(  p_{0}+1\right)
}{\left(  n-1\right)  \left(  p_{0}-r^{-1}\right)  }}\left(  B_{R}^{+}\cup
B_{R}^{n-1}\right)  $. Let%
\[
a=\frac{n\left(  p_{0}+1\right)  }{p_{0}+n-\left(  n-1\right)  p_{0}r},\quad
b=\frac{\left(  p_{0}+1\right)  r}{p_{0}r-1}.
\]
Then $\frac{n}{ra}+\frac{n-1}{b}=\frac{1}{r}$ and%
\[
\frac{r}{\frac{n\left(  p_{0}+1\right)  }{\left(  n-1\right)  p_{0}}}+\frac
{1}{a}=\frac{p_{0}+n}{n\left(  p_{0}+1\right)  }<1.
\]
For $\frac{n\left(  p_{0}+1\right)  }{\left(  n-1\right)  p_{0}}%
<q<\frac{n\left(  p_{0}+1\right)  }{\left(  n-1\right)  \left(  p_{0}%
-r^{-1}\right)  }$, we have $\frac{r}{q}+\frac{1}{a}>\frac{1}{n}$. It follows
from \cite[Proposition 5.2]{HWY} that when $R$ is small enough, $\left.
u_{0}\right\vert _{B_{R/4}^{+}}\in L^{q}\left(  B_{R/4}^{+}\right)  $. This
implies%
\begin{align*}
f_{0}\left(  \xi\right)   & =\int_{B_{R/4}^{+}}P\left(  x,\xi\right)
u_{0}\left(  x\right)  ^{\frac{p_{0}+n}{\left(  n-1\right)  p_{0}}}d\mu\left(
x\right)  +f_{R/4}\left(  \xi\right) \\
& \leq c\left(  M,g,q\right)  \left\vert u_{0}\right\vert _{L^{q}\left(
B_{R/4}^{+}\right)  }^{\frac{p_{0}+n}{\left(  n-1\right)  p_{0}}}%
+f_{R/4}\left(  \xi\right)
\end{align*}
when $q>\frac{n\left(  p_{0}+n\right)  }{\left(  n-1\right)  p_{0}}$. Such a
choice of $q$ is possible since $\frac{n\left(  p_{0}+1\right)  }{\left(
n-1\right)  \left(  p_{0}-r^{-1}\right)  }>\frac{n\left(  p_{0}+n\right)
}{\left(  n-1\right)  p_{0}}$. In particular, we see $\left.  f_{0}\right\vert
_{B_{R/8}}\in L^{\infty}\left(  B_{R/8}\right)  $. Since $\xi_{0}$ is
arbitrary, we see $f_{0}\in L^{\infty}\left(  \Sigma\right)  $ and hence
$u_{0}\in L^{\infty}\left(  M\right)  $. Observing that $f_{0}=T\left(
u_{0}^{\frac{p_{0}+n}{\left(  n-1\right)  p_{0}}}\right)  $, here $T$ is
defined in Section \ref{sec2.1}, it follows from $L^{p}$ theory (\cite[Chapter
9]{GT}) and the Sobolev embedding theorem that $f_{0}\in C^{\alpha}\left(
\Sigma\right)  $ for $0<\alpha<1$. In particular, $f_{0}\left(  \xi\right)
>0$ for any $\xi\in\Sigma$. This implies $u_{0}\in C^{\beta}\left(  M\right)
$ for some $0<\beta<1$ (\cite[Chapter 8]{GT}). It follows from Schauder theory
(\cite[Chapter 6]{GT}) that $f_{0}\in C^{1,\beta}\left(  \Sigma\right)  $.
Iterating this procedure we see $f_{0}\in C^{\infty}\left(  \Sigma\right)  $
and so is $f$.

\begin{case}
\label{case4.2}$\frac{n}{n-1}\leq p_{0}<\infty$.
\end{case}

In this case, we fix a number $r$ such that%
\[
1\leq r\leq p_{0}\text{ and }r\geq\frac{\left(  n-1\right)  p_{0}}{p_{0}+n},
\]
then%
\[
u_{0}\left(  x\right)  ^{1/r}\leq\left[  \int_{B_{R}}P\left(  x,\xi\right)
f_{0}\left(  \xi\right)  ^{p_{0}}dS\left(  \xi\right)  \right]  ^{1/r}%
+u_{R}\left(  x\right)  ^{1/r}.
\]
Hence%
\begin{align*}
& f_{0}\left(  \xi\right) \\
& \leq\int_{B_{R}^{+}}P\left(  x,\xi\right)  u_{0}\left(  x\right)
^{\frac{p_{0}+n}{\left(  n-1\right)  p_{0}}-r^{-1}}\left[  \int_{B_{R}%
}P\left(  x,\zeta\right)  f_{0}\left(  \zeta\right)  ^{p_{0}-r}f_{0}\left(
\zeta\right)  ^{r}dS\left(  \zeta\right)  \right]  ^{1/r}d\mu\left(  x\right)
\\
& +g_{R}\left(  \xi\right) \\
& \leq c\left(  M,g,p,r\right)  \int_{B_{R}^{+}}\frac{x_{n}}{\left(
\left\vert x^{\prime}-\xi\right\vert ^{2}+x_{n}^{2}\right)  ^{n/2}}%
u_{0}\left(  x\right)  ^{\frac{p_{0}+n}{\left(  n-1\right)  p_{0}}-r^{-1}%
}\cdot\\
& \left[  \int_{B_{R}}\frac{x_{n}}{\left(  \left\vert x^{\prime}%
-\zeta\right\vert ^{2}+x_{n}^{2}\right)  ^{n/2}}f_{0}\left(  \zeta\right)
^{p_{0}-r}f_{0}\left(  \zeta\right)  ^{r}d\zeta\right]  ^{1/r}dx+g_{R}\left(
\xi\right)  ,
\end{align*}
here%
\[
g_{R}\left(  \xi\right)  =\int_{B_{R}^{+}}P\left(  x,\xi\right)  u_{0}\left(
x\right)  ^{\frac{p_{0}+n}{\left(  n-1\right)  p_{0}}-r^{-1}}u_{R}\left(
x\right)  ^{1/r}d\mu\left(  x\right)  +f_{R}\left(  \xi\right)  .
\]
We have $g_{R}\in L^{p_{0}+1}\left(  B_{R}\right)  \cap L_{loc}^{q}\left(
B_{R}\right)  $ for any $q<\infty$. Let%
\[
a=\frac{p_{0}+1}{p_{0}-r},\quad b=\frac{n\left(  p_{0}+1\right)  r}{\left(
p_{0}+n\right)  r-\left(  n-1\right)  p_{0}},
\]
then $\frac{n-1}{ra}+\frac{n}{b}=1$, $\frac{r}{p_{0}+1}+\frac{1}{a}%
=\frac{p_{0}}{p_{0}+1}\in\left(  0,1\right)  $. For any $p_{0}+1<q<\infty$, it
follows from \cite[Proposition 5.3]{HWY} that when $R$ is small enough, we
have $f_{0}\in L^{q}\left(  B_{R/4}\right)  $. Since $\xi_{0}$ is arbitrary,
we see $f_{0}\in L^{q}\left(  \Sigma\right)  $ and hence $u_{0}\in
L^{\frac{nq}{\left(  n-1\right)  p_{0}}}\left(  M\right)  $. Using the
equations of $f_{0}$ and $u_{0}$, we see $f_{0}\in L^{\infty}\left(
\Sigma\right)  $ and $u_{0}\in L^{\infty}\left(  M\right)  $. The arguments in
Case \ref{case4.1} tell us $f\in C^{\infty}\left(  \Sigma\right)  $.
\end{proof}

\section{An asymptotic expansion formula of the Poisson kernel\label{sec5}}

Later on we will need more accurate information about the Poisson kernel than
Lemma \ref{lem2.2}. For that purpose we need an asymptotic expansion formula
for this kernel.

Assume $n\geq2$, $\left(  M^{n},g\right)  $ is a smooth compact Riemannian
manifold with boundary $\Sigma=\partial M$, $\delta>0$ is a small number such
that $M_{\delta}=\left\{  x\in M:d\left(  x,\Sigma\right)  \leq\delta\right\}
$ is a tubular neighborhood of $\Sigma$ and $\pi:M_{\delta}\rightarrow\Sigma$
denotes the nearest point projection. For $\xi\in\Sigma$, choose a normal
coordinate for $\Sigma$ at $\xi$, namely $\tau_{1},\cdots,\tau_{n-1}$. Let
$C_{\delta}=\left\{  x\in M_{\delta}:d_{\Sigma}\left(  \pi\left(  x\right)
,\xi\right)  \leq\delta\right\}  $. For $\delta$ small, we have a coordinate
near $\xi$ for $M$ as%
\[
\phi:C_{\delta}\rightarrow\overline{B}_{\delta}^{n-1}\times\left[
0,\delta\right]  :x\mapsto\left(  \tau\left(  \pi\left(  x\right)  \right)
,t\left(  x\right)  \right)  .
\]
It is usually called the Fermi coordinate at $\xi$. We will identify
$C_{\delta}$ with $\overline{B}_{\delta}^{n-1}\times\left[  0,\delta\right]  $
through $\phi$. Denote $r=\left\vert x\right\vert $ and $\theta=\frac
{x}{\left\vert x\right\vert }$.

\begin{theorem}
\label{thm5.1}Under the above set up, we may find $a_{i}\in C^{\infty}\left(
S_{+}^{n-1}\right)  $ with $\left.  a_{i}\right\vert _{\partial S_{+}^{n-1}%
}=0$ for $0\leq i\leq n-1$ and a $\psi\in C^{1,1-\varepsilon}\left(  M\right)
$ (for all $\varepsilon>0$) such that%
\[
P\left(  x,0\right)  =\frac{2}{n\omega_{n}}r^{1-n}\sum_{i=0}^{n-1}r^{i}%
a_{i}\left(  \theta\right)  +\psi\left(  x\right)  \text{ for }x\text{ near
}0\text{.}%
\]
Here $\omega_{n}$ is the volume of the unit ball in $\mathbb{R}^{n}$. Moreover
$a_{0}\left(  \theta\right)  =\theta_{n}=\frac{x_{n}}{\left\vert x\right\vert
}$ and $a_{1}$ is determined by%
\[
\left\{
\begin{array}
[c]{l}%
-\Delta_{S^{n-1}}a_{1}=-H\left(  0\right)  -nH\left(  0\right)  \theta_{n}%
^{2}+2n\left(  n+2\right)  h_{ij}\left(  0\right)  \theta_{i}\theta_{j}%
\theta_{n}^{2}\text{ on }S_{+}^{n-1}\\
\left.  a_{1}\right\vert _{\partial S_{+}^{n-1}}=0
\end{array}
\right.  .
\]
Here $i,j$ runs from $1$ to $n-1$, $h_{ij}$ is the second fundamental form
with respect to inner normal direction and $H$ is the mean curvature.
\end{theorem}

To derive the asymptotic formula, we note that $g=g_{ij}dx_{i}\otimes
dx_{j}+dx_{n}\otimes dx_{n}$. We will use $i,j,k,l$ etc to denote indices
running from $1$ to $n-1$. Calculation shows%
\begin{align}
g_{ij}  & =\delta_{ij}-2h_{ij}\left(  0\right)  x_{n}-\frac{1}{3}\left(
R_{\Sigma}\right)  _{ikjl}\left(  0\right)  x_{k}x_{l}-2h_{ij,k}\left(
0\right)  x_{k}x_{n}\label{eq5.1}\\
& +\left(  -R_{injn}\left(  0\right)  +h_{ik}\left(  0\right)  h_{jk}\left(
0\right)  \right)  x_{n}^{2}+O\left(  r^{3}\right)  ;\nonumber
\end{align}

\begin{align}
g^{ij}  & =\delta_{ij}+2h_{ij}\left(  0\right)  x_{n}+\frac{1}{3}\left(
R_{\Sigma}\right)  _{ikjl}\left(  0\right)  x_{k}x_{l}+2h_{ij,k}\left(
0\right)  x_{k}x_{n}\label{eq5.2}\\
& +\left(  R_{injn}\left(  0\right)  +3h_{ik}\left(  0\right)  h_{jk}\left(
0\right)  \right)  x_{n}^{2}+O\left(  r^{3}\right)  ;\nonumber
\end{align}
and%
\begin{align}
\sqrt{G}  & =1-H\left(  0\right)  x_{n}-\frac{1}{6}\left(  Rc_{\Sigma}\right)
_{ij}\left(  0\right)  x_{i}x_{j}-H_{i}\left(  0\right)  x_{i}x_{n}%
\label{eq5.3}\\
& +\frac{1}{2}\left(  H\left(  0\right)  ^{2}-\left\vert h\left(  0\right)
\right\vert ^{2}-Rc_{nn}\left(  0\right)  \right)  x_{n}^{2}+O\left(
r^{3}\right)  .\nonumber
\end{align}

Note that%
\begin{align*}
\Delta_{g}u  & =\frac{1}{\sqrt{G}}\partial_{i}\left(  g^{ij}\sqrt{G}%
\partial_{j}u\right)  +\frac{1}{\sqrt{G}}\partial_{n}\left(  \sqrt{G}%
\partial_{n}u\right) \\
& =g^{ij}\partial_{ij}u+\partial_{nn}u+\frac{1}{\sqrt{G}}\partial_{i}\left(
g^{ij}\sqrt{G}\right)  \partial_{j}u+\frac{1}{\sqrt{G}}\partial_{n}\left(
\sqrt{G}\right)  \partial_{n}u.
\end{align*}
This and (\ref{eq5.2}), (\ref{eq5.3}) imply that for $\alpha\in\mathbb{R}$ and
$b\in C^{\infty}\left(  S_{+}^{n-1}\right)  $,%
\[
\Delta_{g}\left(  r^{\alpha}b\left(  \theta\right)  \right)  =r^{\alpha-2}
\left[  \Delta_{S^{n-1}}b\left(  \theta\right)  +\alpha\left(  \alpha
+n-2\right)  b\left(  \theta\right)  \right]  +O\left(  r^{\alpha-1}\right)  .
\]
Let $a_{0}\left(  \theta\right)  =\theta_{n}$, then using (\ref{eq5.2}),
(\ref{eq5.3}) we get%
\begin{align*}
& \Delta_{g}\left(  r^{1-n}a_{0}\left(  \theta\right)  \right) \\
& =r^{-n}\left[  -H\left(  0\right)  -nH\left(  0\right)  \theta_{n}%
^{2}+2n\left(  n+2\right)  h_{ij}\left(  0\right)  \theta_{i}\theta_{j}%
\theta_{n}^{2}\right]  +O\left(  r^{-n+1}\right)  .
\end{align*}
Assume for $1\leq k\leq n-1$, we have found $a_{i}\in C^{\infty}\left(
S_{+}^{n-1}\right)  $, vanishing on $\partial S_{+}^{n-1}$ for $0\leq i\leq
k-1$ with%
\[
\Delta_{g}\left(  r^{1-n}\sum_{i=0}^{k-1}a_{i}\left(  \theta\right)
r^{i}\right)  =r^{k-1-n}b_{k-1}\left(  \theta\right)  +O\left(  r^{k-n}%
\right)  ,
\]
then may solve the Dirichlet problem%
\[
\left\{
\begin{array}
[c]{l}%
-\Delta_{S^{n-1}}a_{k}+\left(  k-1\right)  \left(  n-k-1\right)  a_{k}\left(
\theta\right)  =b_{k-1}\left(  \theta\right)  \text{ on }S_{+}^{n-1}\\
\left.  a_{k}\right\vert _{\partial S_{+}^{n-1}}=0
\end{array}
\right.  .
\]
This is possible because $\left(  k-1\right)  \left(  n-k-1\right)  \geq0$.
Then
\[
\Delta_{g}\left(  r^{1-n}\sum_{i=0}^{k}a_{i}\left(  \theta\right)
r^{i}\right)  =O\left(  r^{k-n}\right)  =r^{k-n}b_{k}\left(  \theta\right)
+O\left(  r^{k+1-n}\right)  .
\]
Hence by induction we may find $a_{i}$ for $0\leq i\leq n-1$ such that%
\[
\Delta_{g}\left(  r^{1-n}\sum_{i=0}^{n-1}a_{i}\left(  \theta\right)
r^{i}\right)  =O\left(  r^{-1}\right)  .
\]
Fix a $\eta\in C^{\infty}\left(  \mathbb{R}^{n}\right)  $ such that
$\eta\left(  x\right)  =1$ for $\left\vert x\right\vert \leq\frac{\delta}{4}$
and $\eta\left(  x\right)  =0$ for $\left\vert x\right\vert \geq\frac{\delta
}{2} $. Let $u=\frac{2}{n\omega_{n}}\eta\cdot r^{1-n}\sum_{i=0}^{n-1}%
a_{i}\left(  \theta\right)  r^{i}$, then $\Delta_{g}u=O\left(  r^{-1}\right)
$. We solve%
\[
\left\{
\begin{array}
[c]{l}%
-\Delta_{g}\psi=\Delta_{g}u\text{ on }M\\
\left.  \psi\right\vert _{\partial M}=0
\end{array}
\right.
\]
to find $\psi\in W^{2,n-\varepsilon}\left(  M\right)  $ for all $\varepsilon
>0$. In particular, $\psi\in C^{1,1-\varepsilon}\left(  M\right)  $ for all
$\varepsilon>0$ and the Poisson kernel $P\left(  x,0\right)  =\frac{2}%
{n\omega_{n}}\eta\cdot r^{1-n}\sum_{i=0}^{n-1}a_{i}\left(  \theta\right)
r^{i}+\psi\left(  x\right)  $.

An almost identical argument gives us similar results for the Poisson kernel
of the conformal Laplacian operator.

\begin{proposition}
\label{prop5.2}Under the same set up as in Theorem \ref{thm5.1}. If $n\geq3$
and $\lambda_{1}\left(  L_{g}\right)  >0$, we may find $a_{i}\in C^{\infty
}\left(  S_{+}^{n-1}\right)  $ with $\left.  a_{i}\right\vert _{\partial
S_{+}^{n-1}}=0$ for $0\leq i\leq n-1$ and a $\psi\in C^{1,1-\varepsilon
}\left(  M\right)  $ (for all $\varepsilon>0$) such that%
\[
P_{L}\left(  x,0\right)  =\frac{2}{n\omega_{n}}r^{1-n}\sum_{i=0}^{n-1}%
r^{i}a_{i}\left(  \theta\right)  +\psi\left(  x\right)  \text{ for }x\text{
near }0\text{.}%
\]
Moreover $a_{0}\left(  \theta\right)  =\theta_{n}$ and $a_{1}$ is determined
by%
\[
\left\{
\begin{array}
[c]{l}%
-\Delta_{S^{n-1}}a_{1}=-H\left(  0\right)  -nH\left(  0\right)  \theta_{n}%
^{2}+2n\left(  n+2\right)  h_{ij}\left(  0\right)  \theta_{i}\theta_{j}%
\theta_{n}^{2}\text{ on }S_{+}^{n-1}\\
\left.  a_{1}\right\vert _{\partial S_{+}^{n-1}}=0
\end{array}
\right.  .
\]

\end{proposition}

\section{A criterion for the existence of maximizers\label{sec6}}

We first recall some notations from \cite{HWY}. For $x\in\mathbb{R}_{+}^{n}$,
$\xi\in\mathbb{R}^{n-1}$, the Poisson kernel of the upper half space is%
\[
P\left(  x,\xi\right)  =\frac{2}{n\omega_{n}}\frac{x_{n}}{\left(  \left\vert
x^{\prime}-\xi\right\vert ^{2}+x_{n}^{2}\right)  ^{n/2}}.
\]
Here $x=\left(  x^{\prime},x_{n}\right)  $. For a function $f$ defined on
$\mathbb{R}^{n-1}$, $\left(  Pf\right)  \left(  x\right)  =\int_{\mathbb{R}%
^{n-1}}P\left(  x,\xi\right)  f\left(  \xi\right)  d\xi$. For $1<p<\infty$,
$\left\vert Pf\right\vert _{L^{\frac{np}{n-1}}\left(  \mathbb{R}_{+}%
^{n}\right)  }\leq c_{n,p}\left\vert f\right\vert _{L^{p}\left(
\mathbb{R}^{n-1}\right)  }$, here%
\[
c_{n,p}=\sup\left\{  \left\vert Pf\right\vert _{L^{\frac{np}{n-1}}\left(
\mathbb{R}_{+}^{n}\right)  }:f\in L^{p}\left(  \mathbb{R}^{n-1}\right)
,\left\vert f\right\vert _{L^{p}\left(  \mathbb{R}^{n-1}\right)  }=1\right\}
.
\]

\begin{theorem}
\label{thm6.1}Assume $n\geq2$, $\left(  M^{n},g\right)  $ is a smooth compact
Riemannian manifold with boundary $\Sigma=\partial M$, $1<p<\infty$. Denote%
\[
c_{M,g,p}=\sup\left\{  \left\vert Pf\right\vert _{L^{\frac{np}{n-1}}\left(
M\right)  }:f\in L^{p}\left(  \Sigma\right)  ,\left\vert f\right\vert
_{L^{p}\left(  \Sigma\right)  }=1\right\}  .
\]
Then $c_{M,g,p}\geq c_{n,p}$. Any maximizer of the problem must be smooth and
either strictly positive or strictly negative. Strictly positive maximizers
satisfy the equation%
\[
\int_{M}P\left(  x,\xi\right)  \left(  Pf\right)  \left(  x\right)
^{\frac{np}{n-1}-1}d\mu\left(  x\right)  =c_{M,g,p}^{\frac{np}{n-1}}f\left(
\xi\right)  ^{p-1}.
\]
Moreover if $c_{M,g,p}>c_{n,p}$, then $c_{M,g,p}$ is achieved. Indeed any
maximizing sequence has a convergent subsequence in $L^{p}\left(
\Sigma\right)  $.
\end{theorem}

We use the same notations as in Section \ref{sec2.1}. An ingredient in proving
Theorem \ref{thm6.1} is the following $\varepsilon$-version inequality.

\begin{lemma}
\label{lem6.1}Assume $n\geq2$, $\left(  M^{n},g\right)  $ is a smooth compact
Riemannian manifold with boundary $\Sigma=\partial M$, $1<p<\infty$. Then for
any $\varepsilon>0$ small, there exists a $\delta=\delta\left(
M,g,p,\varepsilon\right)  >0$ such that for every $f\in L^{p}\left(
\Sigma\right)  $,%
\[
\left\vert Pf\right\vert _{L^{\frac{np}{n-1}}\left(  M_{\delta}\right)  }%
\leq\left(  c_{n,p}+\varepsilon\right)  \left\vert f\right\vert _{L^{p}\left(
\Sigma\right)  }.
\]

\end{lemma}

To prove the lemma, we will need the following estimates.

\begin{lemma}
\label{lem6.2}Assume $0\leq\alpha<n-1$, $1<p<\infty$, then%
\[
\left\vert \int_{\Sigma}\frac{f\left(  \xi\right)  }{d\left(  x,\xi\right)
^{\alpha}}dS\left(  \xi\right)  \right\vert _{L^{\frac{np}{\alpha}}\left(
M\right)  }\leq c\left(  M,g,\alpha,p\right)  \left\vert f\right\vert
_{L^{p}\left(  \Sigma\right)  }.
\]

\end{lemma}

\begin{proof}
We may assume $\alpha>0$. For $\varepsilon>0$ small enough, we let
$q_{0}=\frac{n}{\alpha}\left(  1-\varepsilon\right)  $, $q_{1}=1+\frac
{\varepsilon}{p-1}$, then $\frac{1}{p}+\frac{1}{q_{1}}=\frac{q_{0}}{q_{1}%
\cdot\frac{np}{\alpha}}+1$. The needed inequality follows from Lemma
\ref{lem2.3}.
\end{proof}

\begin{corollary}
\label{cor6.1}Assume $\eta\in\operatorname{Lip}\left(  \Sigma\right)  $,
$1<p<\infty$, then%
\[
\left\vert \eta\circ\pi\cdot Pf-P\left(  \eta f\right)  \right\vert
_{L^{\frac{np}{n-2}}\left(  M_{\delta_{0}}\right)  }\leq c\left(
M,g,p\right)  \left\vert \nabla_{\Sigma}\eta\right\vert _{L^{\infty}\left(
\Sigma\right)  }\left\vert f\right\vert _{L^{p}\left(  \Sigma\right)  }.
\]

\end{corollary}

\begin{proof}
It follows from Lemma \ref{lem2.2} that%
\begin{align*}
& \left\vert \eta\left(  \pi\left(  x\right)  \right)  \left(  Pf\right)
\left(  x\right)  -P\left(  \eta f\right)  \left(  x\right)  \right\vert \\
& =\left\vert \int_{\Sigma}\left(  \eta\left(  \pi\left(  x\right)  \right)
-\eta\left(  \xi\right)  \right)  P\left(  x,\xi\right)  f\left(  \xi\right)
dS\left(  \xi\right)  \right\vert \\
& \leq c\left(  M,g\right)  \left\vert \nabla_{\Sigma}\eta\right\vert
_{L^{\infty}\left(  \Sigma\right)  }\int_{\Sigma}\frac{\left\vert f\left(
\xi\right)  \right\vert }{d\left(  x,\xi\right)  ^{n-2}}dS\left(  \xi\right)
.
\end{align*}
Then the conclusion follows from Lemma \ref{lem6.2}.
\end{proof}

\begin{corollary}
\label{cor6.2}Let $K\left(  x,\xi\right)  =\frac{2}{n\omega_{n}}\frac{t\left(
x\right)  }{\left[  t\left(  x\right)  ^{2}+d_{\Sigma}\left(  \pi\left(
x\right)  ,\xi\right)  ^{2}\right]  ^{\frac{n}{2}}}$ for $x\in M_{\delta_{0}}$
and $\xi\in\Sigma$, $\left(  Kf\right)  \left(  x\right)  =\int_{\Sigma
}K\left(  x,\xi\right)  f\left(  \xi\right)  dS\left(  \xi\right)  $,
$1<p<\infty$, then%
\[
\left\vert Pf-Kf\right\vert _{L^{\frac{np}{n-2}}\left(  M_{\delta_{0}}\right)
}\leq c\left(  M,g,p\right)  \left\vert f\right\vert _{L^{p}\left(
\Sigma\right)  }.
\]

\end{corollary}

This follows from Theorem \ref{thm5.1} and Lemma \ref{lem6.2}.

\begin{proof}
[Proof of Lemma \ref{lem6.1}]Without losing of generality we may assume
$f\geq0$. For $\delta_{1}>0$ small, we may find $\eta_{i}\in C^{\infty}\left(
\Sigma,\mathbb{R}\right)  $ for $1\leq i\leq m$ such that $0\leq\eta_{i}\leq
1$, $\sum_{i=1}^{m}\eta_{i}=1$, $\eta_{i}^{1/p}\in C^{\infty}\left(
\Sigma,\mathbb{R}\right)  $ and for each $i$, there exists a point $\xi_{i}%
\in\Sigma$ such that $\eta_{i}\left(  \xi\right)  =0$ for $\xi\in\Sigma$ with
$d_{\Sigma}\left(  \xi,\xi_{i}\right)  \geq\delta_{1}$. For $0<\delta
<\delta_{1}$, we denote%
\[
C_{i,\delta}=\left\{  x\in M_{\delta}:d_{\Sigma}\left(  \pi\left(  x\right)
,\xi_{i}\right)  \leq\delta_{1}\right\}  .
\]
Then%
\begin{align*}
\left\vert Pf\right\vert _{L^{\frac{np}{n-1}}\left(  M_{\delta}\right)  }^{p}
& =\left\vert \left(  Pf\right)  ^{p}\right\vert _{L^{\frac{n}{n-1}}\left(
M_{\delta}\right)  }=\left\vert \sum_{i=1}^{m}\eta_{i}\circ\pi\cdot\left(
Pf\right)  ^{p}\right\vert _{L^{\frac{n}{n-1}}\left(  M_{\delta}\right)  }\\
& \leq\sum_{i=1}^{m}\left\vert \eta_{i}\circ\pi\cdot\left(  Pf\right)
^{p}\right\vert _{L^{\frac{n}{n-1}}\left(  C_{i,\delta}\right)  }=\sum
_{i=1}^{m}\left\vert \eta_{i}^{1/p}\circ\pi\cdot Pf\right\vert _{L^{\frac
{np}{n-1}}\left(  C_{i,\delta}\right)  }^{p}.
\end{align*}
On the other hand, using Corollary \ref{cor6.1} we see%
\begin{align*}
& \left\vert \eta_{i}^{1/p}\circ\pi\cdot Pf\right\vert _{L^{\frac{np}{n-1}%
}\left(  C_{i,\delta}\right)  }\\
& \leq\left\vert P\left(  \eta_{i}^{1/p}f\right)  \right\vert _{L^{\frac
{np}{n-1}}\left(  C_{i,\delta}\right)  }+\left\vert \eta_{i}^{1/p}\circ
\pi\cdot Pf-P\left(  \eta_{i}^{1/p}f\right)  \right\vert _{L^{\frac{np}{n-1}%
}\left(  C_{i,\delta}\right)  }\\
& \leq\left\vert P\left(  \eta_{i}^{1/p}f\right)  \right\vert _{L^{\frac
{np}{n-1}}\left(  C_{i,\delta}\right)  }+\left\vert \eta_{i}^{1/p}\circ
\pi\cdot Pf-P\left(  \eta_{i}^{1/p}f\right)  \right\vert _{L^{\frac{np}{n-2}%
}\left(  C_{i,\delta}\right)  }\left\vert C_{i,\delta}\right\vert ^{\frac
{1}{np}}\\
& \leq\left\vert P\left(  \eta_{i}^{1/p}f\right)  \right\vert _{L^{\frac
{np}{n-1}}\left(  C_{i,\delta}\right)  }+c\left(  M,g,p,\delta_{1}\right)
\delta^{\frac{1}{np}}\left\vert f\right\vert _{L^{p}\left(  \Sigma\right)  }.
\end{align*}
Similarly, by Corollary \ref{cor6.2} we have%
\begin{align*}
& \left\vert P\left(  \eta_{i}^{1/p}f\right)  \right\vert _{L^{\frac{np}{n-1}%
}\left(  C_{i,\delta}\right)  }\\
& \leq\left\vert K\left(  \eta_{i}^{1/p}f\right)  \right\vert _{L^{\frac
{np}{n-1}}\left(  C_{i,\delta}\right)  }+c\left(  M,g,p\right)  \delta
^{\frac{1}{np}}\left\vert f\right\vert _{L^{p}\left(  \Sigma\right)  }\\
& \leq c_{n,p}\left(  1+\varepsilon_{1}\right)  \left\vert \eta_{i}%
^{1/p}f\right\vert _{L^{p}\left(  \Sigma\right)  }+c\left(  M,g,p\right)
\delta^{\frac{1}{np}}\left\vert f\right\vert _{L^{p}\left(  \Sigma\right)  }.
\end{align*}
Here $\varepsilon_{1}=\varepsilon_{1}\left(  M,g,p,\delta_{1}\right)  $ is a
small number which tends to $0$ when $\delta_{1}$ tends to $0$. Hence%
\begin{align*}
\left\vert Pf\right\vert _{L^{\frac{np}{n-1}}\left(  M_{\delta}\right)  }^{p}
& \leq\sum_{i=1}^{m}\left[  c_{n,p}\left(  1+\varepsilon_{1}\right)
\left\vert \eta_{i}^{1/p}f\right\vert _{L^{p}\left(  \Sigma\right)  }+c\left(
M,g,p,\delta_{1}\right)  \delta^{\frac{1}{np}}\left\vert f\right\vert
_{L^{p}\left(  \Sigma\right)  }\right]  ^{p}\\
& \leq\sum_{i=1}^{m}c_{n,p}^{p}\left(  1+2\varepsilon_{1}\right)  ^{p}%
\int_{\Sigma}\eta_{i}f^{p}dS+c\left(  M,g,p,\delta_{1}\right)  \delta
^{1/n}\left\vert f\right\vert _{L^{p}\left(  \Sigma\right)  }^{p}\\
& \leq c_{n,p}^{p}\left(  1+\frac{\varepsilon}{2}\right)  ^{p}\left\vert
f\right\vert _{L^{p}\left(  \Sigma\right)  }^{p}+c\left(  M,g,p,\varepsilon
\right)  \delta^{1/n}\left\vert f\right\vert _{L^{p}\left(  \Sigma\right)
}^{p}\\
& \leq c_{n,p}^{p}\left(  1+\varepsilon\right)  ^{p}\left\vert f\right\vert
_{L^{p}\left(  \Sigma\right)  }^{p}%
\end{align*}
if we first fix $\delta_{1}=\delta_{1}\left(  M,g,p,\varepsilon\right)  $
small enough and then $\delta=\delta\left(  M,g,p,\varepsilon\right)  $ small
enough. This implies $\left\vert Pf\right\vert _{L^{\frac{np}{n-1}}\left(
M_{\delta}\right)  }\leq c_{n,p}\left(  1+\varepsilon\right)  \left\vert
f\right\vert _{L^{p}\left(  \Sigma\right)  }$.
\end{proof}

Next we prove the following concentration compactness lemma (compare with
\cite[lemma 2.1]{L} and \cite[proposition 3.1]{HWY}).

\begin{proposition}
[Concentration compactness lemma]\label{prop6.1}Assume $n\geq2$, $\left(
M^{n},g\right)  $ is a smooth compact Riemannian manifold with boundary
$\Sigma=\partial M$, $1<p<\infty$, $f_{i}\in L^{p}\left(  \Sigma\right)  $
such that $f_{i}\rightharpoonup f$ in $L^{p}\left(  \Sigma\right)  $. After
passing to a subsequence assume%
\[
\left\vert f_{i}\right\vert ^{p}dS\rightharpoonup\sigma\text{ in }%
\mathcal{M}\left(  \Sigma\right)  ,\quad\left\vert Pf_{i}\right\vert
^{\frac{np}{n-1}}d\mu\rightharpoonup\nu\text{ in }\mathcal{M}\left(  M\right)
.
\]
Here $\mathcal{M}\left(  \Sigma\right)  $ is the space of all Radon measures
on $\Sigma$. Then we have

\begin{itemize}
\item $\left.  \nu\right\vert _{M\backslash\Sigma}=\left\vert Pf\right\vert
^{\frac{np}{n-1}}d\mu$. Moreover for every Borel set $E\subset\Sigma$,
$\nu\left(  E\right)  ^{\frac{n-1}{np}}\leq c_{n,p}\sigma\left(  E\right)
^{\frac{1}{p}}$.

\item There exists a countable set of points $\zeta_{j}\in\Sigma$ such that
$\nu=\left\vert Pf\right\vert ^{\frac{np}{n-1}}d\mu+\sum_{j}\nu_{j}%
\delta_{\zeta_{j}}$, $\sigma\geq\left\vert f\right\vert ^{p}dS+\sum_{j}%
\sigma_{j}\delta_{\zeta_{j}}$, here $\sigma_{j}=\sigma\left(  \left\{
\zeta_{j}\right\}  \right)  $ and $\nu_{j}^{\frac{n-1}{np}}\leq c_{n,p}%
\sigma_{j}^{\frac{1}{p}}$.
\end{itemize}
\end{proposition}

\begin{proof}
Without losing of generality we may assume $\left\vert f_{i}\right\vert
_{L^{p}\left(  \Sigma\right)  }\leq1$. Since $\left\vert \left(
Pf_{i}\right)  \left(  x\right)  \right\vert \leq c\left(  M,g,p\right)
t\left(  x\right)  ^{-\frac{n-1}{p}}$ for $x\in M\backslash\Sigma$, it follows
from the elliptic estimates of harmonic functions that $Pf_{i}\rightarrow Pf$
in $C_{loc}^{\infty}\left(  M\backslash\Sigma\right)  $. In particular,
$\left.  \nu\right\vert _{M\backslash\Sigma}=\left\vert Pf\right\vert
^{\frac{np}{n-1}}d\mu$. For $\varepsilon>0$ small, it follows from Lemma
\ref{lem6.1} and Corollary \ref{cor6.1} that for $\varphi\in C^{\infty}\left(
\Sigma\right)  $ and $\delta>0$ small enough,%
\begin{align*}
& \left\vert \varphi\circ\pi\cdot Pf_{i}\right\vert _{L^{\frac{np}{n-1}%
}\left(  M_{\delta}\right)  }\\
& \leq\left\vert P\left(  \varphi f_{i}\right)  \right\vert _{L^{\frac
{np}{n-1}}\left(  M_{\delta}\right)  }+\left\vert \varphi\circ\pi\cdot
Pf_{i}-P\left(  \varphi f_{i}\right)  \right\vert _{L^{\frac{np}{n-1}}\left(
M_{\delta}\right)  }\\
& \leq\left(  c_{n,p}+\varepsilon\right)  \left\vert \varphi f_{i}\right\vert
_{L^{p}\left(  \Sigma\right)  }+c\left(  M,g,p\right)  \delta^{\frac{1}{np}%
}\left\vert \nabla_{\Sigma}\varphi\right\vert _{L^{\infty}\left(
\Sigma\right)  }.
\end{align*}
Let $i\rightarrow\infty$ we see%
\[
\left(  \int_{\Sigma}\left\vert \varphi\right\vert ^{\frac{np}{n-1}}%
d\nu\right)  ^{\frac{n-1}{np}}\leq\left(  c_{n,p}+\varepsilon\right)  \left(
\int_{\Sigma}\left\vert \varphi\right\vert ^{p}d\sigma\right)  ^{\frac{1}{p}%
}+c\left(  M,g,p\right)  \delta^{\frac{1}{np}}\left\vert \nabla_{\Sigma
}\varphi\right\vert _{L^{\infty}\left(  \Sigma\right)  }.
\]
Let $\delta\rightarrow0^{+}$ and then $\varepsilon\rightarrow0^{+}$, we get%
\[
\left(  \int_{\Sigma}\left\vert \varphi\right\vert ^{\frac{np}{n-1}}%
d\nu\right)  ^{\frac{n-1}{np}}\leq c_{n,p}\left(  \int_{\Sigma}\left\vert
\varphi\right\vert ^{p}d\sigma\right)  ^{\frac{1}{p}}.
\]
A limit process shows for every nonnegative Borel function $h$ on $\Sigma$,%
\[
\left(  \int_{\Sigma}h^{\frac{np}{n-1}}d\nu\right)  ^{\frac{n-1}{np}}\leq
c_{n,p}\left(  \int_{\Sigma}h^{p}d\sigma\right)  ^{\frac{1}{p}}.
\]
In particular, for every Borel set $E\subset\Sigma$, $\nu\left(  E\right)
^{\frac{n-1}{np}}\leq c_{n,p}\sigma\left(  E\right)  ^{\frac{1}{p}}$. Based on
this inequality we may proceed as in the proof of \cite[proposition 3.1]{HWY}
to get the second conclusion.
\end{proof}

Now we are ready to derive Theorem \ref{thm6.1}.

\begin{proof}
[Proof of Theorem \ref{thm6.1}]First we want to show $c_{M,g,p}\geq c_{n,p}$
is always true. To see this we may fix a point $\xi_{0}\in\Sigma$, choose a
normal coordinate for $\Sigma$ at $\xi_{0}$, namely $\tau_{1},\cdots
,\tau_{n-1}$. For $\delta>0$ small, we denote $C_{\delta}=\left\{  x\in
M_{\delta}:d_{\Sigma}\left(  \pi\left(  x\right)  ,\xi_{0}\right)  \leq
\delta\right\}  $, then we have a natural coordinate near $\xi_{0}$ for $M$ as%
\[
\phi:C_{\delta}\rightarrow\overline{B}_{\delta}^{n-1}\times\left[
0,\delta\right]  :x\mapsto\left(  \tau\left(  \pi\left(  x\right)  \right)
,t\left(  x\right)  \right)  .
\]
We will identify $C_{\delta}$ with $\overline{B}_{\delta}^{n-1}\times\left[
0,\delta\right]  $ through $\phi$. On $C_{\delta}$ we have the Euclidean
metric $g_{0}=\sum_{i=1}^{n}dx_{i}\otimes dx_{i}$. If $\overline{f}\in
L^{p}\left(  \Sigma\right)  \backslash\left\{  0\right\}  $ and $\overline{f}$
vanishes outside $\overline{B}_{\delta}^{n-1}$, then it follows from Corollary
\ref{cor6.2} that%
\[
\left\vert K\overline{f}\right\vert _{L^{\frac{np}{n-1}}\left(  C_{\delta
},g\right)  }\leq\left\vert P\overline{f}\right\vert _{L^{\frac{np}{n-1}%
}\left(  C_{\delta},g\right)  }+c\left(  M,g,p\right)  \delta^{\frac{1}{p}%
}\left\vert \overline{f}\right\vert _{L^{p}\left(  \Sigma\right)  }.
\]
Let $f\left(  \xi\right)  =\overline{f}\left(  \xi\right)  $ for $\left\vert
\xi\right\vert \leq\delta$ and $f\left(  \xi\right)  =0$ for $\left\vert
\xi\right\vert >\delta,\xi\in\mathbb{R}^{n-1}$, and $u$ be the harmonic
extension of $f$ to $\mathbb{R}_{+}^{n}$, then%
\begin{align*}
\left\vert u\right\vert _{L^{\frac{np}{n-1}}\left(  C_{\delta},g_{0}\right)
}  & \leq\left(  1+\varepsilon_{1}\right)  \left\vert K\overline{f}\right\vert
_{L^{\frac{np}{n-1}}\left(  C_{\delta},g\right)  }\\
& \leq\left(  1+\varepsilon_{1}\right)  \left\vert P\overline{f}\right\vert
_{L^{\frac{np}{n-1}}\left(  C_{\delta},g\right)  }+c\left(  M,g,p\right)
\delta^{\frac{1}{p}}\left\vert \overline{f}\right\vert _{L^{p}\left(
\Sigma\right)  }.
\end{align*}
Here $\varepsilon_{1}=\varepsilon_{1}\left(  M,g,p,\delta\right)  $ and
$\varepsilon_{1}\rightarrow0^{+}$ as $\delta\rightarrow0^{+}$. Hence%
\begin{align*}
c_{M,g,p}  & \geq\frac{\left\vert P\overline{f}\right\vert _{L^{\frac{np}%
{n-1}}\left(  M\right)  }}{\left\vert \overline{f}\right\vert _{L^{p}\left(
\Sigma\right)  }}\geq\frac{\left\vert P\overline{f}\right\vert _{L^{\frac
{np}{n-1}}\left(  C_{\delta},g\right)  }}{\left\vert \overline{f}\right\vert
_{L^{p}\left(  B_{\delta}^{n-1},g\right)  }}\\
& \geq\frac{1}{\left(  1+\varepsilon_{1}\right)  ^{2}}\frac{\left\vert
u\right\vert _{L^{\frac{np}{n-1}}\left(  C_{\delta},g_{0}\right)  }%
}{\left\vert f\right\vert _{L^{p}\left(  B_{\delta}^{n-1},g_{0}\right)  }%
}-c\left(  M,g,p\right)  \delta^{\frac{1}{p}}.
\end{align*}
Assume $f\in L^{p}\left(  \mathbb{R}^{n-1}\right)  \backslash\left\{
0\right\}  $ and $f=0$ outside a ball, $u$ is the harmonic extension of $f$ to
$\mathbb{R}_{+}^{n}$, then for $\varepsilon>0$ small enough, we write
$f_{\varepsilon}\left(  \xi\right)  =\varepsilon^{-\frac{n-1}{p}}f\left(
\frac{\xi}{\varepsilon}\right)  $ and $u_{\varepsilon}\left(  x\right)
=\varepsilon^{-\frac{n-1}{p}}u\left(  \frac{x}{\varepsilon}\right)  $. Let
$\overline{f}=f_{\varepsilon}$ on $B_{\delta}^{n-1}$ and $0$ on $\Sigma
\backslash B_{\delta}^{n-1}$, then we get%
\[
c_{M,g,p}\geq\frac{1}{\left(  1+\varepsilon_{1}\right)  ^{2}}\frac{\left\vert
u_{\varepsilon}\right\vert _{L^{\frac{np}{n-1}}\left(  C_{\delta}%
,g_{0}\right)  }}{\left\vert f_{\varepsilon}\right\vert _{L^{p}\left(
B_{\delta}^{n-1},g_{0}\right)  }}-c\left(  M,g,p\right)  \delta^{\frac{1}{p}}.
\]
Let $\varepsilon\rightarrow0^{+}$ then $\delta\rightarrow0^{+}$, we see%
\[
c_{M,g,p}\geq\frac{\left\vert u\right\vert _{L^{\frac{np}{n-1}}\left(
\mathbb{R}_{+}^{n}\right)  }}{\left\vert f\right\vert _{L^{p}\left(
\mathbb{R}^{n-1}\right)  }}.
\]
By approximation we know the inequality remains true for all $f\in
L^{p}\left(  \mathbb{R}^{n-1}\right)  \backslash\left\{  0\right\}  $ and this
implies $c_{M,g,p}\geq c_{n,p}$.

If $f$ is a maximizer, then it is clear that $f$ will be either nonnegative or
nonpositive. Assume $f\geq0$, then it satisfies the Euler-Lagrange equation%
\[
\int_{M}P\left(  x,\xi\right)  \left(  Pf\right)  \left(  x\right)
^{\frac{np}{n-1}-1}d\mu\left(  x\right)  =c_{M,g,p}^{\frac{np}{n-1}}f\left(
\xi\right)  ^{p-1}.
\]
It follows from Proposition \ref{prop4.1} that $f$ must be smooth and hence it
is strictly positive.

Assume $c_{M,g,p}>c_{n,p}$. Let $f_{i}\in L^{p}\left(  \Sigma\right)  $ be a
sequence of functions with $\left\vert f_{i}\right\vert _{L^{p}\left(
\Sigma\right)  }=1$ and $\left\vert Pf_{i}\right\vert _{L^{\frac{np}{n-1}%
}\left(  M\right)  }\rightarrow c_{M,g,p}$. After passing to a subsequence we
may assume $f_{i}\rightharpoonup f$ in $L^{p}\left(  \Sigma\right)  $,
$\left\vert f_{i}\right\vert ^{p}dS\rightharpoonup\sigma$ in $\mathcal{M}%
\left(  \Sigma\right)  $ and $\left\vert Pf_{i}\right\vert ^{\frac{np}{n-1}%
}d\mu\rightharpoonup\nu$ in $\mathcal{M}\left(  M\right)  $. It follows from
Proposition \ref{prop6.1} that we may find a countable set of points
$\zeta_{j}\in\Sigma$ such that $\nu=\left\vert Pf\right\vert ^{\frac{np}{n-1}%
}d\mu+\sum_{j}\nu_{j}\delta_{\zeta_{j}}$ and $\sigma\geq\left\vert
f\right\vert ^{p}dS+\sum_{j}\sigma_{j}\delta_{\zeta_{j}}$. Here $\sigma
_{j}=\sigma\left(  \left\{  \zeta_{j}\right\}  \right)  $ and $\nu_{j}%
^{\frac{n-1}{n}}\leq c_{n,p}^{p}\sigma_{j}$. In particular $1=\sigma\left(
\Sigma\right)  \geq\left\vert f\right\vert _{L^{p}\left(  \Sigma\right)  }%
^{p}+\sum_{j}\sigma_{j}$. We claim $\nu_{j}=0$ for all $j$. If this is not the
case, then%
\[
c_{M,g,p}^{\frac{np}{n-1}}=\nu\left(  M\right)  =\left\vert Pf\right\vert
_{L^{\frac{np}{n-1}}\left(  M\right)  }^{\frac{np}{n-1}}+\sum_{j}\nu_{j}\leq
c_{M,g,p}^{\frac{np}{n-1}}\left\vert f\right\vert _{L^{p}\left(
\Sigma\right)  }^{\frac{np}{n-1}}+\sum_{j}\nu_{j}.
\]
Hence%
\begin{align*}
c_{M,g,p}^{p}  & \leq c_{M,g,p}^{p}\left\vert f\right\vert _{L^{p}\left(
\Sigma\right)  }^{p}+\sum_{j}\nu_{j}^{\frac{n-1}{n}}\leq c_{M,g,p}%
^{p}\left\vert f\right\vert _{L^{p}\left(  \Sigma\right)  }^{p}+c_{n,p}%
^{p}\sum_{j}\sigma_{j}\\
& <c_{M,g,p}^{p}\left\vert f\right\vert _{L^{p}\left(  \Sigma\right)  }%
^{p}+c_{M,g,p}^{p}\sum_{j}\sigma_{j}.
\end{align*}
This implies $1<\left\vert f\right\vert _{L^{p}\left(  \Sigma\right)  }%
^{p}+\sum_{j}\sigma_{j}$, a contradiction. Since $\nu_{j}=0$ for all $j$, we
see $\left\vert Pf\right\vert _{L^{\frac{np}{n-1}}\left(  M\right)
}=c_{M,g,p}$. Hence $\left\vert f\right\vert _{L^{p}\left(  \Sigma\right)
}\geq1$. This implies $f_{i}\rightarrow f$ in $L^{p}\left(  \Sigma\right)  $.
That is every maximizing sequence has a convergent subsequence in
$L^{p}\left(  \Sigma\right)  $ and $c_{M,g,p}$ is achieved.
\end{proof}

\section{Proof of the Theorem \ref{thm1.1}\label{sec7}}

In this section we finish the proof of Theorem \ref{thm1.1}. Without losing of
generality we may assume $R=0$. It follows from Theorem \ref{thm3.1}, Theorem
\ref{thm6.1} and \cite[thm1.1]{HWY} that%
\[
\Theta_{M,g}=c_{M,g,\frac{2\left(  n-1\right)  }{n-2}}^{\frac{2}{n-2}}\geq
c_{n,\frac{2\left(  n-1\right)  }{n-2}}^{\frac{2}{n-2}}=n^{-\frac{1}{n-1}%
}\omega_{n}^{-\frac{1}{n\left(  n-1\right)  }}=\Theta_{\overline{B}%
_{1},g_{\mathbb{R}^{n}}}.
\]
On the other hand, if $\Theta_{M,g}>\Theta_{\overline{B}_{1},g_{\mathbb{R}%
^{n}}}$, then $c_{M,g,\frac{2\left(  n-1\right)  }{n-2}}>c_{n,\frac{2\left(
n-1\right)  }{n-2}}$. It follows from Theorem \ref{thm6.1} that we may find a
$f\in C^{\infty}\left(  \Sigma\right)  $ with $f>0$ such that $\left\vert
f\right\vert _{L^{\frac{2\left(  n-1\right)  }{n-2}}\left(  \Sigma\right)
}=1$ and $c_{M,g,\frac{2\left(  n-1\right)  }{n-2}}=\left\vert Pf\right\vert
_{L^{\frac{2n}{n-2}}\left(  M\right)  }$. Let $\widetilde{g}=\left(
Pf\right)  ^{\frac{4}{n-2}}g$, then clearly $\widetilde{R}=0$ and $I\left(
M,\widetilde{g}\right)  =\Theta_{M,g}$.

\end{document}